\documentclass[12pt]{amsart}
\usepackage{graphicx}
\usepackage{amssymb}
\usepackage[shortlabels]{enumitem}
\usepackage[margin=1.25in]{geometry}
\newtheorem{thm}{Theorem}[section]

\theoremstyle{definition}

\theoremstyle{remark}
\newtheorem{rem}[thm]{Remark}

\begin{document}
\title[]{Bernstein Theorems for Nonlinear Geometric PDEs}
\author[]{Connor Mooney}
\maketitle
\tableofcontents

\newpage
\section{Introduction and Acknowledgments}
In 1915, Bernstein proved a beautiful theorem which says that the only entire minimal graphs in $\mathbb{R}^3$ are planes. The problem of classifying global solutions has since become a driving force in the study of nonlinear geometric PDEs. 

\vspace{3mm}

In this article we revisit the Bernstein problem for several geometric PDEs including the minimal surface, Monge-Amp\`{e}re, and special Lagrangian equations. We also discuss the minimal surface system where appropriate. We first explore equations in two variables. We then discuss fully nonlinear equations in higher dimensions. Third, we discuss rigidity results for the minimal surface equation in higher dimensions. Finally, we discuss examples of non-flat entire minimal graphs in high dimensions. Our exposition includes a construction of the celebrated Bombieri-De Giorgi-Giusti example, using the methods that were recently introduced to solve the Bernstein problem for anisotropic minimal hypersurfaces \cite{MY}. We feel that this approach highlights intuition for readers learning the subject. In the last section we survey recent results and state some open questions about variants of the Bernstein problem.

\vspace{3mm}

The article is based on a lecture series given by the author for the inaugural European Doctorate School of Differential Geometry, held in Granada in June 2024. The author is grateful to Jos\'{e} Espinar, Jos\'{e} G\'{a}lvez, Francisco Torralbo, and Magdalena Rodr\'{i}guez for organizing the event. This work was supported by a Simons Fellowship, a Sloan Research Fellowship, and NSF CAREER Grant DMS-2143668.

\newpage
\section{Louville Theorems for Uniformly Elliptic PDEs}\label{UnifElliptic}

The basic Liouville theorem says that a nonnegative harmonic function $u$ on $\mathbb{R}^n$ is constant. If $u$ vanishes at a point, the result follows immediately from the strong maximum principle. The idea is to ``quantify" this approach. After adding a constant we may assume that $\inf_{\mathbb{R}^n} u = 0$. For any $\epsilon > 0$, we may assume after a translation that $u(0) < \epsilon$. The mean value property gives a dimensional constant $C$ such that 
$$\int_{B_2(0)} u \leq Cu(0) \text{ and } u(y) \leq C\int_{B_1(y)} u$$ 
for any $y \in B_1$, hence $u < C^2\epsilon$ in $B_1$. The same reasoning applied to $u(r \cdot)$ for any $r > 0$ gives $u < C^2\epsilon$ globally, and since $\epsilon$ was arbitrary we are done.

The same result holds if instead $Lu = 0$, where $L = \text{div}(A(x)\nabla \cdot)$ (divergence form) or $L = \text{tr}(A(x)D^2\cdot)$ (non divergence form), and $A(x)$ is a symmetric matrix with eigenvalues in $[\lambda,\,\lambda^{-1}]$ for some fixed $\lambda > 0$ independent of $x$. Since no further regularity of $A$ is assumed, such operators cannot be treated as perturbations of the Laplace operator (like in Schauder theory). Nonetheless, such operators enjoy versions of the mean value property:
\begin{equation}\label{SubHarnack}
w \geq 0 \text{ and } Lw \geq 0 \text{ in } B_1 \Rightarrow \sup_{B_{1/2}}w^p \leq C(n,\,\lambda,\,p)\int_{B_1} w^p \text{ for any } p > 0,
\end{equation}
\begin{equation}\label{SuperHarnack}
w \geq 0 \text{ and } Lw \leq 0 \text{ in } B_1 \Rightarrow \int_{B_{1/2}} w^{\delta} \leq C(n,\,\lambda) w^{\delta}(0) \text{ for some } \delta(n,\,\lambda) > 0.
\end{equation}
These are deep results of De Giorgi \cite{DG} and Nash \cite{Na} in the case that $L$ has divergence form, and of Krylov-Safonov \cite{KS} in the case $L$ is in non divergence form. Proofs can be found in the canonical references \cite{GT} (both forms) and \cite{CC2} (non divergence form). The Liouville theorem (global solutions to $Lu = 0$ that are bounded from one side are constants) follows in the same way as for harmonic functions.

As a consequence we get rigidity theorems for global solutions to some nonlinear PDEs. We first discuss the quasilinear case. Assume that $H$ is a smooth function on $\mathbb{R}^n$ and $D^2H > 0$. Assume further that $u$ has bounded gradient, and solves
\begin{equation}\label{Quasilinear}
\text{div}(\nabla H(\nabla u)) = 0
\end{equation}
on $\mathbb{R}^n$. Equivalently, $u$ is a minimizer of the energy $\int H(\nabla u)$. Then the derivatives $u_k$ of $u$ solve a divergence-form equation of the type discussed above, with $A = D^2H(\nabla u)$ and $\lambda > 0$ depending on $H$ and $\sup_{\mathbb{R}^n}|\nabla u|$. Thus, $u$ is linear.

The fully nonlinear case is trickier. Assume that $F$ is a smooth function of the $n \times n$ symmetric matrices $M$ such that $DF = (F_{ij}) = (\partial_{M_{ij}}F)$ is everywhere positive definite. Assume further that $u$ has bounded Hessian and solves
\begin{equation}\label{FullyNonlinearEqn}
F(D^2u) = 0
\end{equation}
on $\mathbb{R}^n$. We would like to conclude that $u$ is a quadratic polynomial. This turns out to be false, at least in dimensions $n \geq 5$ (see \cite{NTV}). (It remains open in dimensions three and four, and the result is true in two dimensions by work of Nirenberg \cite{Nir}, see also Remark \ref{2D}). It is true under the additional assumption that $\{F = 0\}$ is a convex hypersurface, which is satisfied in important applications. This result is due to Evans \cite{E} and Krylov \cite{Kr}. We now give a proof using the idea in \cite{CSil}.

 The first and (pure) second derivatives $u_k,\,u_{kk}$ of $u$ solve
$$Lu_k = DF \cdot D^2u_k = 0, \quad Lu_{kk} = -D^2F(D^2u_k,\,D^2u_k),$$
where $L$ is a non divergence equation of the type discussed above with $A = DF(D^2u)$ and $\lambda > 0$ depending on $F$ and $\sup_{\mathbb{R}^n}|D^2u|$. Because $\{F = 0\}$ is a convex hypersurface, and since $D^2u_k$ is tangent to $\{F = 0\}$ by the once-differentiated equation, $Lu_{kk}$ has a sign. We may assume it is nonnegative after possibly replacing $u$ by $-u$ and $F$ by $-F(-\cdot)$. However, the ellipticity of the equation philosophically means that $u_{kk}$ is a negative combination of other second derivatives, which suggests that $u_{kk}$ behaves like a solution to $Lw \leq 0$ as well. We make this reasoning rigorous now.

We may assume that $D^2u(0) = 0$ after subtracting a quadratic polynomial from $u$ and translating $F$. For a symmetric matrix $M$, we let $G_+(M)$, resp. $G_-(M)$ denote the sum of its positive, resp. negative eigenvalues. Note that, by the fundamental theorem of calculus,
$$0 = \int_0^1 \frac{d}{dt} F(tD^2u)\,dt = \left(\int_0^1 F_{ij}(tD^2u)\,dt\right)u_{ij}.$$
This implies that $C^{-1}G_+(D^2u) \leq -G_-(D^2u) \leq CG_+(D^2u)$. Here and below $C$ will denote a large positive constant depending only on $n,\,\lambda$, and may change from line to line. It thus suffices to prove that $G_+(D^2u) \equiv 0$.

If not, then up to multiplying $u$ by a constant and making a Lipschitz rescaling of $F$ we may assume that $\sup_{\mathbb{R}^n} G_+(D^2u) = 1$. After making a quadratic rescaling of $u$, we may assume that $\sup_{B_1} G_+(D^2u) > 1-\epsilon$, for $\epsilon$ as small as we like to be chosen. Thus there is a point in $B_1$ and a subspace $V$ of $\mathbb{R}^n$ where $\Delta_Vu > 1-\epsilon$. Here $\Delta_V$ denotes the Laplace operator on $V$, i.e. $\sum_{i = 1}^m u_{z_iz_i}$, where $\{z_i\}_{i = 1}^m$ form an orthonormal basis for $V$. Noting that $L(1-\Delta_Vu) \leq 0$ and applying (\ref{SuperHarnack}), appropriately translated and rescaled, we get that, in $B_1$, $1-\epsilon^{1/2} \leq \Delta_Vu \leq G_+(D^2u) \leq 1$
away from a set $S$ of measure at most $C\epsilon^{\delta/2}$.

\begin{figure}
 \centering
    \includegraphics[scale=0.7, trim={0mm 125mm 0mm 20mm}, clip]{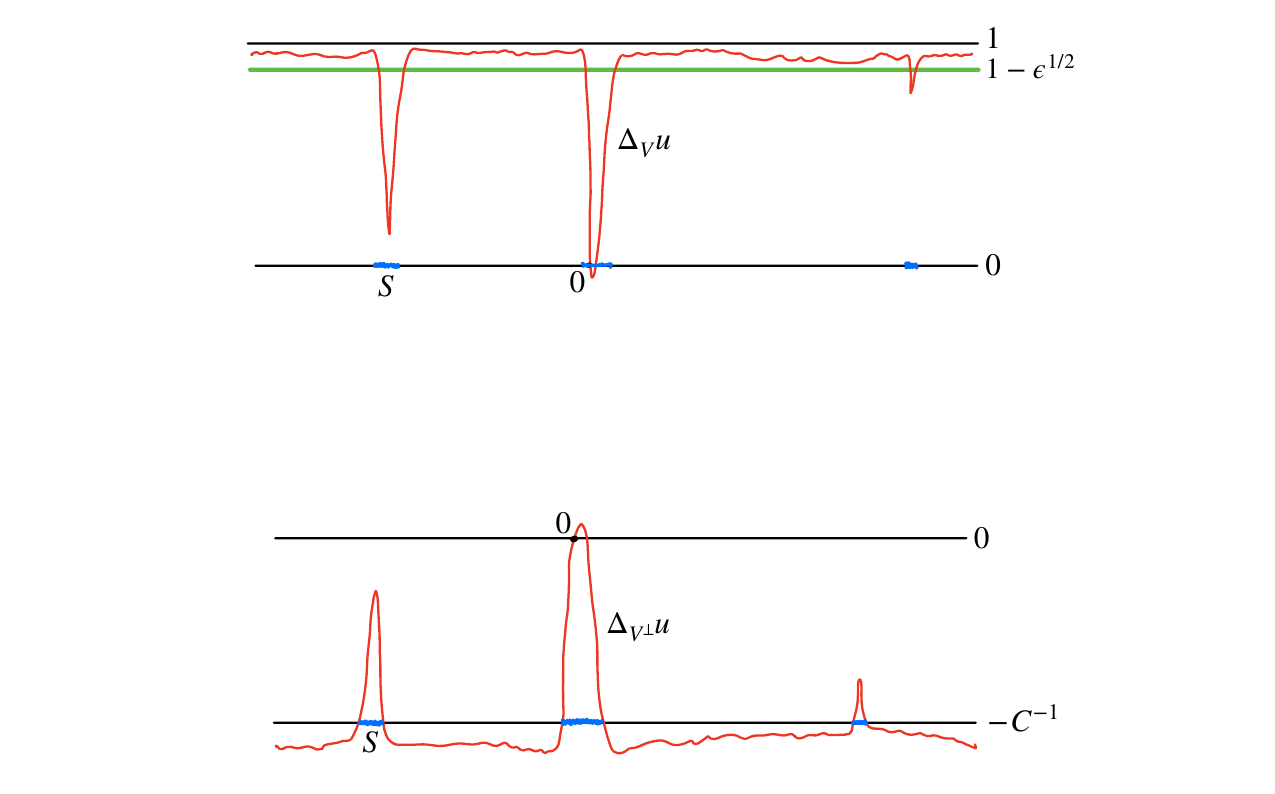}
\caption{}
\label{EKFig}
\end{figure}

Now let $V^{\perp}$ denote the subspace orthogonal to $V$. From 
$$\Delta u = G_+(D^2u) + G_-(D^2u) = \Delta_V u + \Delta_{V^{\perp}} u$$ 
we see that $\Delta_{V^{\perp}} u$ is within $\epsilon^{1/2}$ of $G_-(D^2u)$ away from $S$. Recalling that $G_- \leq -C^{-1}G_+$, we conclude for $\epsilon$ small that $\Delta_{V^{\perp}} u \leq -C^{-1}$ away from $S$. Applying (\ref{SubHarnack}) to $w = (\Delta_{V^{\perp}} u + C^{-1})_+$ with $p = 1$ and noting that $w = 0$ away from $S$, we get
$$C^{-1} = w(0) \leq C\int_{B_1} w \leq C|S| \leq C\epsilon^{\delta/2},$$
a contradiction for $\epsilon(n,\,\lambda) > 0$ small. See Figure \ref{EKFig}.

\begin{rem}\label{2D}
In two variables we have stronger results. For example, if $L$ has divergence form, $Lw \leq 0$ on $\mathbb{R}^2$, and $w \geq 0$, then $w$ is constant. The proof goes as follows. Let $\psi$ be a compactly supported Lipschitz function on the plane. For any fixed $k > 0$ use $(k-w)_+ \psi^2$ as a test function and apply Cauchy-Schwarz to get the Caccioppoli inequality
$$\int_{\mathbb{R}^2} |\nabla (k-w)_+|^2\psi^2 \leq C(\lambda,\,k)\int_{\mathbb{R}^2} |\nabla \psi|^2.$$
Taking $\psi = 1$ in $B_1$ and $\max\{1-\log|x|/\log(R),\,0\}$ outside $B_1$, the right hand side is bounded above by $C(\log R)^{-1}$, and taking $R$ to $\infty$ gives that $(k-w)_+$ is a constant for any $k > 0$. The result follows. The point is that, in two dimensions, one pays no energy to cut off using $\log$. The proof is in the same spirit as the proofs of the Bernstein theorem for minimal surfaces based on the stability inequality, see Section \ref{MSEHigherD}. It is a good exercise to prove that there are non-constant positive superharmonic functions in $\mathbb{R}^3$, so this result is false in higher dimensions. The result is also false for non-divergence form equations in the plane, as can be seen from the example $w = (1+|x|^2)^{-1}$, which has Hessian eigenvalues $\lambda_1 \geq \lambda_2$ with $\lambda_2 < 0$ and $-3\lambda_2 > |\lambda_1|$. 

In the non divergence setting, the correct analogue is: If $L$ is a non divergence form linear, uniformly elliptic operator, then solutions to $Lw = 0$ with bounded gradient on the plane are linear. This implies the Liouville theorem for (\ref{FullyNonlinearEqn}) when $n = 2$ without any additional hypothesis on $\{F = 0\}$ needed (apply the result in the previous sentence to the first derivatives of $u$). To prove the result, it suffices after rotation to prove that $w_1$ is constant. After dividing by $A_{22}$ the equation can be written as a divergence-form, uniformly elliptic equation $\text{div}(\tilde{A}\nabla w_1) = 0$, where $\tilde{A}_{11} = A_{11}/A_{22},\, \tilde{A}_{12} = 2A_{12}/A_{22}, \, \tilde{A}_{21} = 0,\, \tilde{A}_{22} = 1$. The Liouville theorem for divergence-form equations implies that $w_1$ is constant. (The coefficients aren't symmetric, but the relevant thing is that $\tilde{A} + \tilde{A}^T$ has eigenvalues bounded between fixed positive constants). When $n \geq 4$ the corresponding result is false, e.g. the one-homogeneous function $(|z_1|^2 - |z_2|^2)/|z|$ on $\mathbb{C}^2 \cong \mathbb{R}^4$ solves a non-divergence form linear, uniformly elliptic equation. On the other hand, the only one-homogeneous solutions to such equations in dimension $n = 3$ are linear functions (see \cite{A}, \cite{HNY}; there are solutions that are homogeneous of degree $\alpha \in (0,\,1)$ due to Safonov, see \cite{Saf}). Morally, one-homogeneity makes the problem two-dimensional. The question whether bounded gradient implies linear for global solutions to non divergence form linear, uniformly elliptic equations on $\mathbb{R}^3$ remains open.

The aforementioned result in \cite{A}, \cite{HNY} implies that the only two-homogeneous functions $u$ that solve fully nonlinear, uniformly elliptic equations of the form (\ref{FullyNonlinearEqn}) in $\mathbb{R}^3$ are quadratic polynomials (apply the result to the first derivatives of $u$). Interestingly, the same is true with $\mathbb{R}^3$ replaced by $\mathbb{R}^4$, provided $u$ is assumed to be analytic away from the origin \cite{NV1}. The latter result is dimensionally optimal \cite{NTV}.
\end{rem}

\newpage
\section{Equations in Two Variables}\label{2DNonlinear}
\subsection{Monge-Amp\`{e}re Equation}
Let $u$ be a smooth solution to
$$\det D^2u = 1$$
in $\mathbb{R}^2$. We will show that $u$ is a quadratic polynomial, no growth condition needed. This was proven by J\"{o}rgens \cite{J}.

\begin{figure}
 \centering
    \includegraphics[scale=0.7, trim={0mm 140mm 0mm 20mm}, clip]{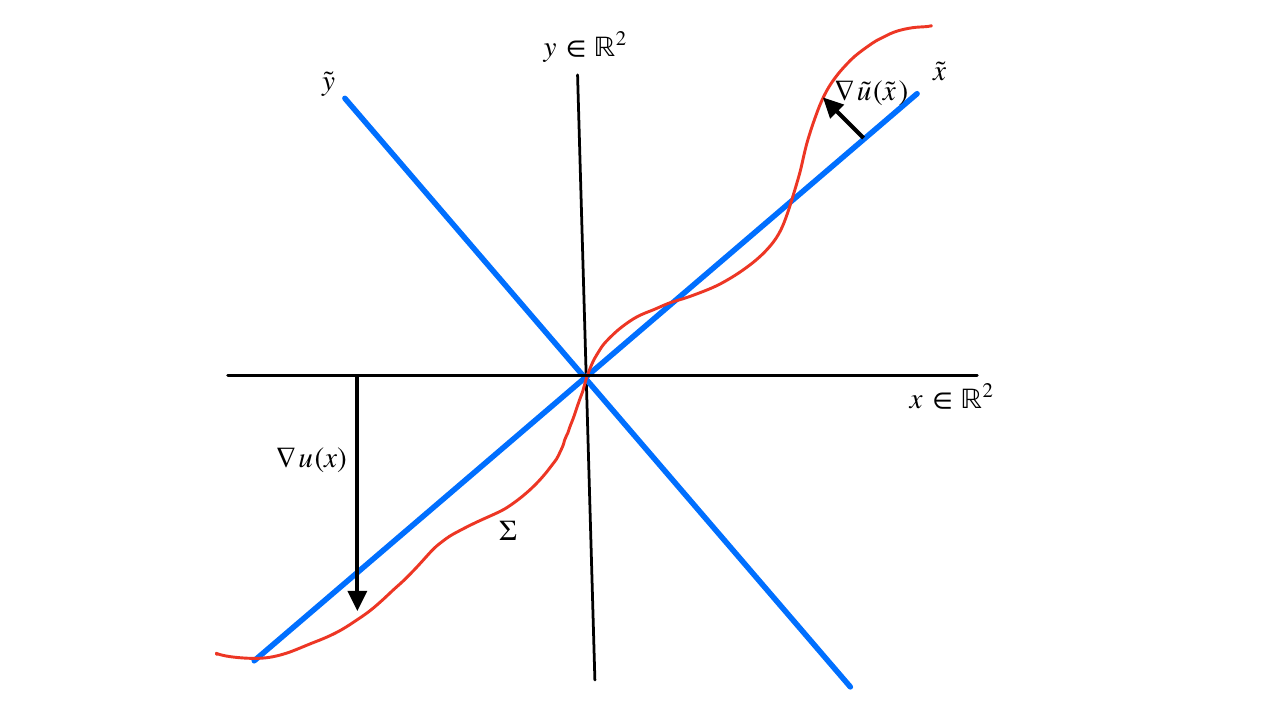}
\caption{}
\label{LegLewy}
\end{figure}

Up to replacing $u$ by $-u$ we may assume that $u$ is convex. We exploit the Legendre-Lewy transform. That is, we study the equations solved by potentials of the rotated gradient graph. We write the gradient graph $\Sigma := \{(x,\,\nabla u(x))\} \subset \mathbb{R}^4$ in terms of the $SU(2)\, \pi/4$-rotated coordinates
$$\tilde{x} = \frac{1}{\sqrt{2}} x + \frac{1}{\sqrt{2}}y, \quad \tilde{y} = -\frac{1}{\sqrt{2}}x + \frac{1}{\sqrt{2}}y, \quad x,\,y \in \mathbb{R}^2 \quad \text{(see Figure \ref{LegLewy}).}$$
That is, $\Sigma = \{(\tilde{x}, G(\tilde{x}))\}$, where
$$G\left(x/\sqrt{2} + \nabla u(x)/\sqrt{2}\right) := -x/\sqrt{2} + \nabla u(x)/\sqrt{2}.$$
Note that $x + \nabla u(x) = \nabla (|x|^2/2 + u)$ is a diffeomorphism of $\mathbb{R}^2$ by the convexity of $u$, hence $\Sigma$ is a global graph over $\{\tilde{y} = 0\}$. Differentiating, it is easy to see that $DG$ is symmetric, hence $G = \nabla \tilde{u}$ for some function $\tilde{u}$, and moreover
$$D^2\tilde{u} = (-I + D^2u)(I + D^2u)^{-1}.$$
The eigenvalues of $D^2\tilde{u}$ are thus $\frac{-1+\lambda_1}{1+\lambda_1} = \frac{-\lambda_2 + 1}{\lambda_2 + 1}, \, \frac{-1+\lambda_2}{1+\lambda_2}$, where $\lambda_i$ are the eigenvalues of $D^2u$ (in the first equality we used the equation $\lambda_1\lambda_2 = 1$). We conclude that $\tilde{u}$ is harmonic with bounded Hessian (in fact, $-I < D^2\tilde{u} < I$), hence it is a quadratic polynomial. In particular, $\Sigma$ is a $2$-plane, as desired.

\begin{rem}\label{LLT}
The Legendre-Lewy transform works more generally as follows. Assume that $D^2u \geq -\cot(\theta + \delta)I = \tan(\theta + \delta - \pi/2)I$ in $\mathbb{R}^n$ with $\theta,\,\delta$ positive and $\theta + \delta < \pi$. Then write the gradient graph $\Sigma$ of $u$ in terms of the rotated coordinates
$$\tilde{x} = \cos\theta \,x + \sin\theta \,y, \quad \tilde{y} = -\sin\theta \,x + \cos\theta \,y,\, x,\,y \in \mathbb{R}^n.$$
The lower bound on $D^2u$ guarantees that $\Sigma$ is an entire graph over $\{\tilde{y} = 0\}$, since $\cot(\theta) |x|^2/2 + u$ is uniformly convex. In the same way as above, we conclude that $\Sigma$ is the gradient graph of a new potential $\tilde{u}$, whose Hessian eigenvalues $\tilde{\lambda_k}$ are related to the Hessian eigenvalues $\lambda_k$ of $u$ by
\begin{equation}\label{NewEigenvals}
\tilde{\lambda_k} = \frac{-\sin\theta + \cos\theta \lambda_k}{\cos\theta + \sin\theta \lambda_k} \Longleftrightarrow \tan^{-1}\tilde{\lambda_k} = \tan^{-1}\lambda_k - \theta.
\end{equation}
Since $-\pi/2 + \theta + \delta \leq \tan^{-1}\lambda_k < \pi/2,$ we conclude that 
$$-\pi/2 + \delta \leq \tan^{-1}\tilde{\lambda_k} < \pi/2-\theta,$$
that is, $\tilde{u}$ has bounded Hessian. We note that in the case $\theta = \frac{\pi}{2}$, $D^2\tilde{u} = -(D^2u)^{-1}$ and $\tilde{u}$ is the negative Legendre transform of $u$.

Using this tool we can prove the Liouville theorem for a family of equations in two variables that includes the Monge-Amp\`{e}re equation as a special case. If $u$ is a smooth function on $\mathbb{R}^2$ such that
$$\sum_{k \leq 2} \tan^{-1}\lambda_k  = \mu \neq 0,$$
then $u$ is a quadratic polynomial. The above equation is called the special Lagrangian equation, and it says that the gradient graph of $u$ is a minimal surface (of codimension two). The case $\mu = \pi/2$ is the Monge-Amp\`{e}re equation. To prove the Liouville theorem, we may assume after possibly replacing $u$ by $-u$ that $\mu > 0$. Applying the Legendre-Lewy rotation with $\delta = \theta = \mu/2$ gives us a harmonic function with bounded Hessian, hence, flat gradient graph, as desired. We generalize this result to higher dimensions in Section \ref{FullyNonlinear}.
\end{rem}

\subsection{Minimal Surface Equation}
To begin the discussion we find it convenient to consider the case of general dimension and codimension. A smooth map 
$$u = (u^1,\,...,\,u^m) : \mathbb{R}^n \rightarrow \mathbb{R}^m$$ 
solves the minimal surface system if it is a critical point of the area functional
$$A(u) := \int \sqrt{\det g(Du)}\,dx,$$
where $g(p) = I_{n \times n} + p^Tp,\, p \in \mathbb{R}^{m \times n}$. Smooth critical points of integrals of the form $E(u) := \int F(Du)$, with $F$ smooth on $\mathbb{R}^{m \times n}$, satisfy that $E(u + \epsilon \varphi) = E(u) + o(\epsilon)$, i.e.
\begin{equation}\label{Variation}
\int F_{p^{\alpha}_i}(Du) \varphi^{\alpha}_i \,dx = 0,
\end{equation}
for all compactly supported smooth maps $\varphi$. Here $F_{p^{\alpha}_i}$ denotes the derivative of $F$ in the $p^{\alpha}_i$ direction of $\mathbb{R}^{m \times n}$. This is typically written as the system
$$\partial_{i}(F_{p^{\alpha}_i}(Du)) = 0,\, \alpha = 1,\,...,\,m,$$
known as the outer variation system, since it reflects the criticality of $u$ under perturbations of its values. In the case of the area functional $F = \sqrt{\det g}$ this is
\begin{equation}\label{Outer}
\partial_i(\sqrt{\det g}g^{ij}\partial_ju^{\alpha}) = 0,\, \alpha = 1,\,...,\,m.
\end{equation}

\begin{figure}
 \centering
    \includegraphics[scale=0.7, trim={0mm 200mm 0mm 20mm}, clip]{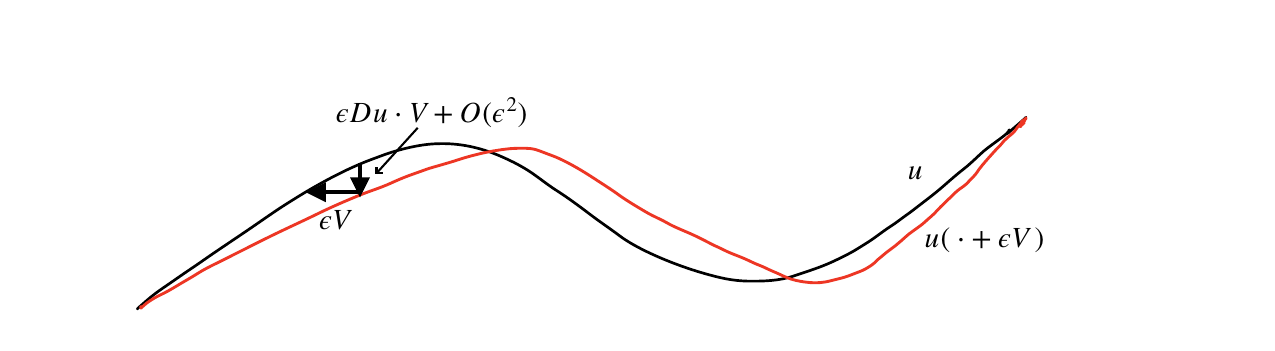}
\caption{}
\label{InnerOuter}
\end{figure}

It is sometimes convenient to exploit the criticality of $u$ under deformations of the domain instead: $E(u(x + \epsilon V)) = E(u) + o(\epsilon)$ for compactly supported $V$. Since $u(x + \epsilon V) = u(x) + \epsilon Du \cdot V + O(\epsilon^2)$, the domain variation $V$ is equivalent to the outer variation $Du \cdot V$ (see Figure \ref{InnerOuter}). Taking $\varphi = Du \cdot V$ in (\ref{Variation}) and integrating by parts leads to the system
$$\partial_i(F_{p^{\alpha}_i}(Du)u^{\alpha}_j - F(Du)\delta_{ij}) = 0, \quad j = 1,\,...,\,n.$$
This is called the inner variation system. In the case of the area functional $F = \sqrt{\det g}$, this is equivalent to
\begin{equation}\label{Inner}
\partial_i(\sqrt{\det g}g^{ij}) = 0,\, j = 1,\,...,\,n.
\end{equation}
We leave this calculation as an exercise. From the derivation it is clear that (\ref{Outer}) implies (\ref{Inner}). To conclude the general discussion we note that (\ref{Inner}) and (\ref{Outer}) imply that
\begin{equation}\label{MSS}
g^{ij}u^{\alpha}_{ij} = 0,\, \alpha = 1,\,...,\,m.
\end{equation}
The system (\ref{MSS}) is in fact equivalent to (\ref{Outer}). This can be seen by using the fact that, for any map $v$, $\partial_i(\sqrt{\det g}g^{ij}\partial_j(x,\,v(x)))$ is normal to the graph $\{(x,\,v(x))\}$ of $v$.

\begin{rem}
The above calculations use that $u$ has two derivatives. It is not clear whether Lipschitz maps that solve (\ref{Outer}) in the sense of distributions also solve (\ref{Inner}). In the codimension one case $m = 1$, Lipschitz solutions to (\ref{Outer}) are smooth by well-known elliptic regularity theory (see e.g. \cite{GT}) so there is no issue. Lawson and Osserman conjectured that Lipschitz solutions to (\ref{Outer}) also solve (\ref{Inner}) in higher codimension in their seminal work \cite{LO}. This conjecture was confirmed in the case $n = 2$, $m$ arbitrary only recently \cite{HMT}, and it remains open when $n \geq 3,\, m \geq 2$.
\end{rem}

We now specialize to the case that $n = 2$ and $m$ is arbitrary. The inner variation system reads $\text{div}(\sqrt{\det g}g^{-1}) = 0$, with the divergence taken row-wise. Since divergence-free vector fields in the plane are rotations of gradients, we have
$$\sqrt{\det g}g^{-1} = DwJ, \quad J := \left(\begin{array}{cc}
0 & -1 \\
1 & 0
\end{array}\right),$$
for some map $w$. It follows that $D(Jw) = JDw = -J\sqrt{\det g}g^{-1}J$. The symmetry of the RHS implies that $Jw = \nabla \Phi$ for some potential $\Phi$. We conclude that
$$D^2\Phi = JDw = -J(\sqrt{\det g}g^{-1})J.$$
The matrix on the RHS is positive with determinant one, so $\Phi$ is a solution to the Monge-Amp\`{e}re equation. The J\"{o}rgens theorem implies that $\sqrt{\det g}g^{-1}$ is constant.

The above considerations worked in arbitrary codimension. In the codimension one case $m = 1$, the eigenvalues of $\sqrt{\det g}g^{-1}$ are $(1 + |\nabla u|^2)^{\pm 1/2}$, so it follows that $|\nabla u|$ is constant. Proving that $\nabla u$ is constant is a simple exercise from here, yielding Bernstein's theorem \cite{Be}: a solution to the minimal surface equation on the whole plane must be linear.

\begin{rem}
The above approach is due to Nitsche \cite{Ni}. The original proof of Bernstein \cite{Be} is more involved. It is based on two results. The first is that $w := \tan^{-1}u_k$ is harmonic on the graph of $u$, $k = 1,\,2$. This can be seen by noting that the unit normal $\nu = (-\nabla u,\,1)/\sqrt{1+|\nabla u|^2}$ is anti-conformal, as its differential is the second fundamental form of the graph, which is symmetric and trace-free. Thus, $w$ can be viewed as the phase of a holomorphic map. As a consequence, $w$ is a bounded function on $\mathbb{R}^2$ that solves an elliptic equation of the form $a_{ij}w_{ij} = 0$, where $(a_{ij}) > 0$. The second remarkable result is that any such function must be constant (no quantitative information about $a_{ij}$ needed). This is established using a topological argument.
\end{rem}

\subsection{Minimal Surface System}
We conclude the discussion of equations in two variables with the minimal surface system, i.e. $n = 2$ and $m \geq 2$. We proved above that $\sqrt{\det g}g^{-1}$ is constant. After a rotation of the plane we may assume that $\sqrt{\det g}g^{-1}$ is diagonal with positive entries $\lambda,\,\lambda^{-1}$. Thus $u$ solves
$$\lambda u_{11} + \lambda^{-1}u_{22} = 0, \quad u_1 \cdot u_2 = 0, \quad \frac{1+|u_2|^2}{1+|u_1|^2} = \lambda^2.$$
Letting $\tilde{u}(x_1,\,x_2) := u(\sqrt{\lambda} x_1,\, x_2/\sqrt{\lambda}),$
this becomes
$$\Delta \tilde{u} = 0, \quad \tilde{u}_1 \cdot \tilde{u}_2 = 0, \quad |\tilde{u}_2|^2 - |\tilde{u}_1|^2 = \lambda - \lambda^{-1} := \Lambda.$$
In terms of $h := \tilde{u}_2 + i\tilde{u}_1$, the system takes the convenient form
\begin{equation}\label{HolMSS}
h \text{ holomorphic}, \quad \sum_{\alpha = 1}^m (h^{\alpha})^2 = \Lambda.
\end{equation}
In the case $m = 2$ and $\Lambda \neq 0$ we can proceed further. The second equation in (\ref{HolMSS}) can be written
$$(h^1 + ih^2)(h^1 - ih^2) = \Lambda \neq 0.$$
Thus, $h^1 + ih^2 = e^{H},\, h^1 - ih^2 = \Lambda e^{-H}$ for a holomorphic $H$, that is,
\begin{equation}\label{Cod2}
h^1 = \frac{1}{2}(e^{H} + \Lambda e^{-H}), \quad h^2 = \frac{1}{2i}(e^{H} - \Lambda e^{-H}).
\end{equation}

The preceding discussion allows us to classify all global solutions $u$ to the minimal surface system in the case $n = m = 2$. The first possibility is that $u$ is holomorphic or anti-holomorphic, corresponding to the case $\lambda = 1$. The remaining solutions can be obtained as follows. For any $\lambda \in (0,\,\infty) \backslash \{1\}$ and holomorphic $H$, define $h$ as in (\ref{Cod2}) with $\Lambda := \lambda - \lambda^{-1}$. Then, integrate the curl-free vector fields $(\text{Im}(h^{\alpha}),\,\text{Re}(h^{\alpha}))$ to get $\tilde{u}^{\alpha}$ up to constants. Next let $v := \tilde{u}(x_1/\sqrt{\lambda},\,\sqrt{\lambda} x_2)$. Finally, let $u$ be $v$ composed with any rotation. This characterization of entire graphical minimal surfaces of dimension two in $\mathbb{R}^4$ can be found in the survey of Osserman \cite{Oss}.
 
 \newpage
 \section{Fully Nonlinear Equations in Higher Dimensions}\label{FullyNonlinear}
  
 \subsection{Special Lagrangian Equation}
 The special Lagrangian equation is
 \begin{equation}\label{SLAG}
 F_{\Theta}(D^2u) := \sum_{i = 1}^n \tan^{-1}(\lambda_i) - \Theta = 0,
 \end{equation}
 where $u$ is defined on a domain in $\mathbb{R}^n$ and $\lambda_i$ denote the Hessian eigenvalues of $u$. There is a rigidity theorem for global solutions when $|\Theta|$ is sufficiently large, due to Yuan \cite{Y2}. We prove it here.
 
 The point is that the set $\{F_{\Theta} = 0\}$ is convex if and only if $|\Theta| \geq (n-2)\pi/2$. Here is the calculation. First, it suffices to consider the case $\Theta \geq (n-2)\pi/2$, since $\{F_{\pm\Theta} = 0\}$ are reflections through the origin of one another. Let $H(x) = \sum_{i = 1}^n \tan^{-1}(x_i)$ on $\mathbb{R}^n$. It suffices to prove that $\{H = \Theta\}$ is convex (see e.g. \cite{CNS}, Section 3). Let $z$ be tangent to $\{H = \Theta\}$, that is,
$$\sum_{i = 1}^n z_i(1+x_i^2)^{-1} = 0.$$
We aim to show that 
\begin{equation}\label{LevSetCvx}
-\frac{1}{2}D^2H(z,\,z) = \sum_{i = 1}^n x_iz_i^2(1+x_i^2)^{-2} \geq 0.
\end{equation}
If all $x_i \geq 0$ this is obvious, so assume that $x_1 \geq ... \geq x_n$ and $x_n < 0$. Then by the equation $H(x) = \Theta \geq (n-2)\pi/2 + \delta$ with $\delta \geq 0$ we have $x_{n-1} > 0$ and $\delta < \pi/2$. The only term in (\ref{LevSetCvx}) that we worry about is the last, which we can estimate by the first order condition and Cauchy-Schwarz:
$$z_n^2(1+x_n^2)^{-2} \leq \left(\sum_{i < n} x_iz_i^2(1+x_i^2)^{-2}\right)\left(\sum_{i < n} x_i^{-1}\right).$$
Using this inequality in (\ref{LevSetCvx}) reduces the problem to proving that
$$\sum_{i = 1}^n x_i^{-1} \leq 0.$$
To that end note that 
$$\sum_{i = 1}^n \tan^{-1}(x_i^{-1}) = \sum_{i < n} (\pi/2 -\tan^{-1}(x_i)) -\pi/2 - \tan^{-1}x_n = -\delta,$$
so $\tan\left(\sum_{i = 1}^n \tan^{-1}(x_i^{-1})\right) \leq 0$ and $0 < \sum_{i < n} \tan^{-1}(x_i^{-1}) < \pi/2-\delta$. The sum formula for tangent then implies that
$$\tan\left(\sum_{i < n} \tan^{-1}(x_i^{-1})\right) + x_n^{-1} \leq 0,$$
and the result follows from the superadditivity of $\tan$.

A rigidity theorem (global solutions are quadratic polynomials) follows for solutions to $F_{\Theta}(D^2u) = 0$, provided $|\Theta| = (n-2)\pi/2 + \mu$ for some $\mu> 0$. Indeed, up to replacing $u$ by $-u$, we may assume that $\Theta = (n-2)\pi/2 + \mu$. Perform the Legendre-Lewy transform with $\theta = \delta = \mu/n$ (see Remark \ref{LLT}). This gives a function $\tilde{u}$ with bounded Hessian on $\mathbb{R}^n$ that solves $F_{(n-2)\pi/2}(D^2\tilde{u}) = 0$. Since the level set $\{F_{(n-2)\pi/2} = 0\}$ is convex, the general Liouville theorem for fully nonlinear PDEs from Section \ref{UnifElliptic} applies. We conclude that the gradient graph of $\tilde{u}$ (hence that of $u$) is flat, as desired.

\begin{rem}
It is important that $\mu > 0$. For example, the case $n = 2,\, \mu = 0$ corresponds to the Laplace equation, which admits many entire solutions. In the case $n = 3,\, \mu = \pi/2$ there are global solutions with exponential growth \cite{W}.
\end{rem}

\begin{rem}
There are rigidity results for entire solutions to $F_{\Theta}(D^2u) = 0$ under other conditions as well. These results use that (\ref{SLAG}) is the potential equation for minimal graphs that are half the dimension of the ambient space, i.e. that $\nabla u$ solves the minimal surface system. For example, if $u$ is convex, or if $u$ is semi-convex (Hessian bounded below) and $n \leq 4$, then it is a quadratic polynomial \cite{Y1}. The idea is to perform a Legendre-Lewy rotation by an appropriate angle to get a new potential $\tilde{u}$ with bounded Hessian, and apply results for Lipschitz solutions to the minimal surface system (see Section \ref{MSSBernstein} below for a discussion of the relevant results). When $u$ is convex, rotation by $\pi/4$ gives a potential with $-I \leq D^2\tilde{u} < I$, and $\nabla \tilde{u}$ satisfies the area-decreasing condition. (It is not the strict one described in Section \ref{MSSBernstein}, so more work is required, but this is the idea). When $u$ is semi-convex and $n \leq 4$, rotation by an arbitrary small angle gives the desired potential. Indeed, in the case $n \leq 3$, Lipschitz entire solutions to the minimal surface system are linear. When $n = 4$ the general theory of the minimal surface system doesn't suffice (see e.g. the Lawson-Osserman example in Section \ref{LOEx}). However, the Lagrangian structure can be exploited. The monotonicity formula allows one to reduce to the case that $\tilde{u}$ is homogeneous of degree two and analytic away from the origin, and rigidity follows from the main result in \cite{NV1}. Rigidity for entire semi-convex solutions to (\ref{SLAG}) in general dimension $n$ would follow from the non-existence of non-flat graphical special Lagrangian cones of dimension $n$, but that remains unknown when $n \geq 5$.
\end{rem}

\begin{rem}
The technique of Legendre-Lewy rotation has also been effective for studying other fully nonlinear PDEs. For example, it was used to prove the rigidity of convex, resp. semiconvex global solutions to $\sigma_2(D^2u) = 1$ on $\mathbb{R}^n$, see \cite{CY}, resp. \cite{SY}.
\end{rem}

\subsection{Monge-Amp\`{e}re equation}
A well-known result is that global convex solutions in $\mathbb{R}^n$ to
$$\det D^2u = 1$$
are quadratic polynomials. This is due to J\"{o}rgens when $n = 2$ (proof above), Calabi \cite{Cal} for $n \leq 5$, and Pogorelov \cite{Pog} in all dimensions. We outline the proof below. 

When $n \geq 3$ the Monge-Amp\`{e}re equation is no longer an instance of the special Lagrangian equation, so new ideas are needed. The equation can be written $\sum_{k = 1}^n \log \lambda_k = 0$. The concavity of $\log$ implies that the branch of $\{\det = 1\}$ in the positive symmetric matrices is convex, so it suffices by the general discussion in Section \ref{UnifElliptic} to get a global Hessian bound. 

The key ingredients to do so are the following. 

First, the equation is affine invariant: $u(A\cdot)$ also solves the equation, provided $|\det A| = 1$.

\begin{figure}
 \centering
    \includegraphics[scale=0.7, trim={0mm 110mm 0mm 30mm}, clip]{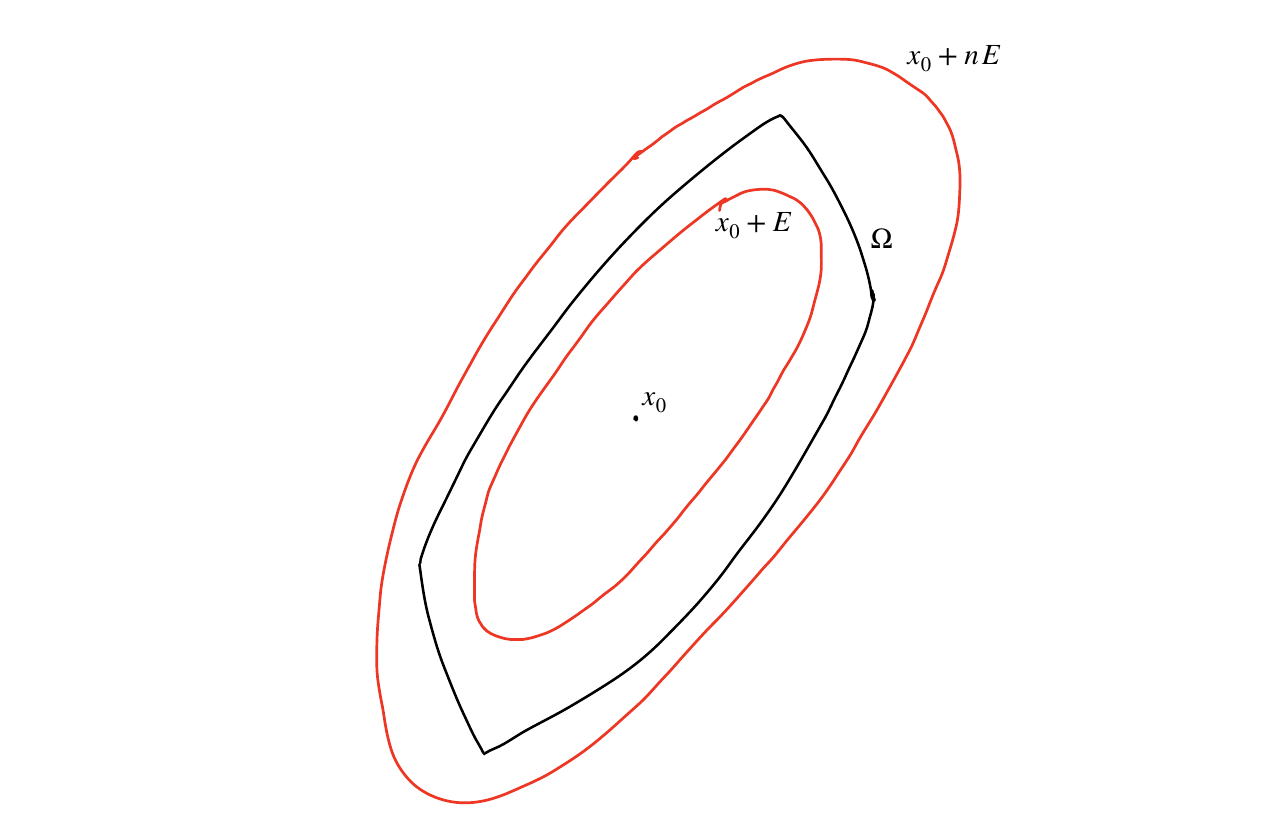}
\caption{}
\label{JohnsLemma}
\end{figure}

The second ingredient, which allows us to exploit affine invariance, is John's lemma: for any bounded open convex set $\Omega \subset \mathbb{R}^n$, there exists an ellipsoid $E$ centered at $0$ and a point $x_0 \in \Omega$ such that
$$x_0 + E \subset \Omega \subset x_0 + nE$$
(see Figure \ref{JohnsLemma}). The ellipsoid $x_0 + E$ is the ellipsoid of maximal volume contained in $\Omega$, and the factor $n$ is sharp on simplices. 

The last ingredient is Pogorelov's interior $C^2$ estimate: if $w$ is a convex solution to $\det D^2w = 1$ in $\{w < 0\}$ and $B_{\kappa} \subset \{w < 0\} \subset B_{\kappa^{-1}}$ for some $\kappa > 0$, then for $\delta > 0$ we have
$$\|D^2w\|_{L^{\infty}(\{w < -\delta\})} \leq C(n,\,\kappa,\,\delta).$$
To prove this Pogorelov considers the point where the quantity $Q = w_{kk}|w|e^{w_k^2/2}$ attains its maximum ($k$ is any direction), and uses the information that at this point, $\nabla Q = 0$ and $LQ \leq 0$, where $L$ is the linearized operator $w^{ij}\partial_i\partial_j$. The first two derivatives of the equation itself are used to simplify expressions.

We continue with the proof of rigidity of global solutions. Up to subtracting a linear function and performing an affine transformation, we may assume that $u(0) = |\nabla u(0)| = 0,\, D^2u(0) = I$. It suffices to prove that there exists $c(n) > 0$ such that
\begin{equation}\label{MAContainment}
B_c \subset \{u < 1\} \subset B_{c^{-1}}.
\end{equation}
Indeed, assume that (\ref{MAContainment}) is true. Then the same holds with $u$ replaced by $u_R := R^{-2}u(R \cdot)$, since $u_R$ satisfies the same equation and conditions. Applying the Pogorelov estimate to $w = u_R - 1$ with $\delta = 1/2$ gives an upper bound for $|D^2u|$ in $\{u < R^2/2\}$ independent of $R$, so the desired Hessian bound follows after sending $R$ to $\infty$. 

\begin{figure}
 \centering
    \includegraphics[scale=0.7, trim={0mm 75mm 0mm 20mm}, clip]{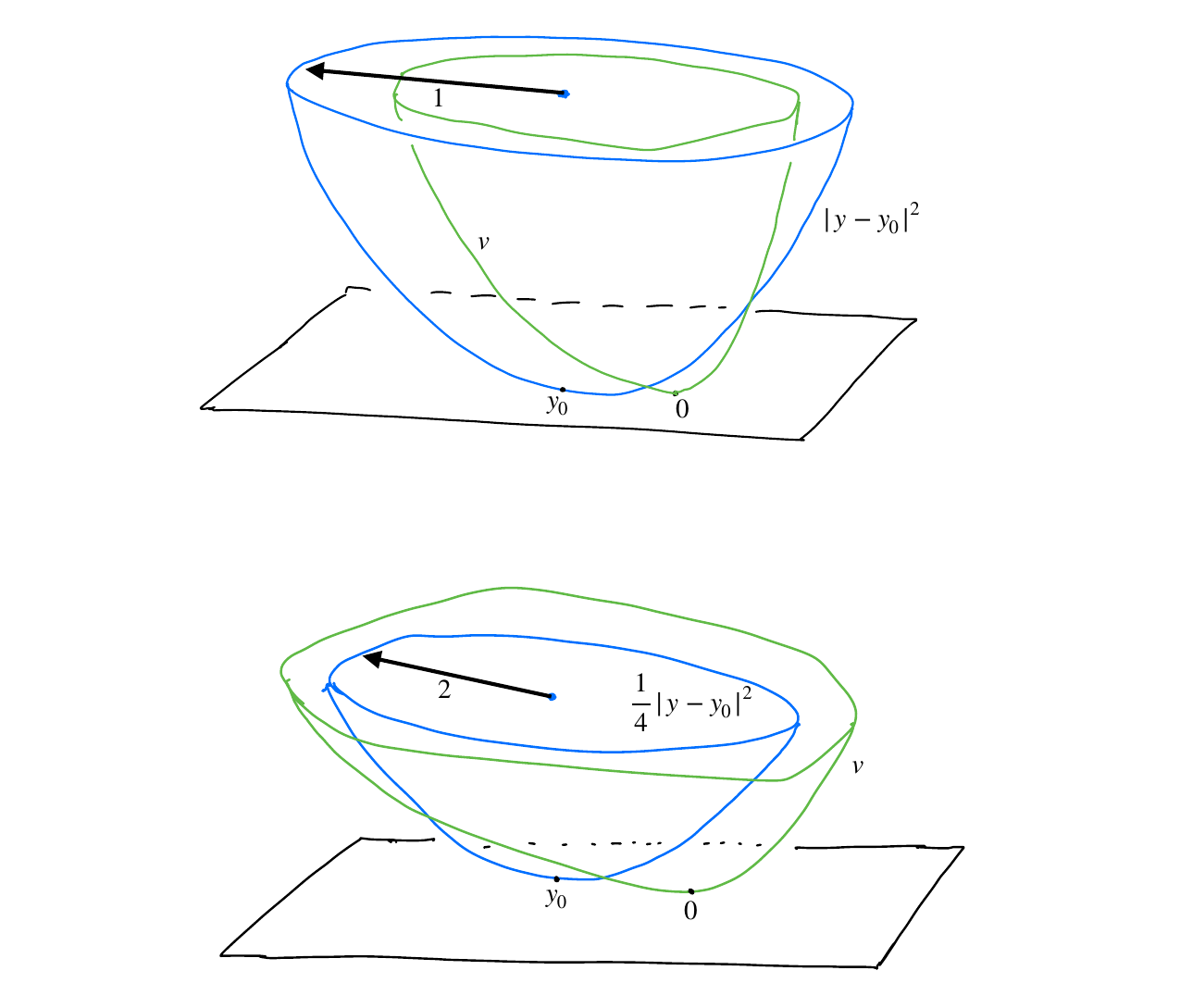}
\caption{}
\label{MAMaxPr}
\end{figure}

We proceed with the proof of (\ref{MAContainment}). By John's Lemma, there is a volume-preserving affine transformation $A$, a number $k > 0$, and a point $y_0$ such that for $v = u(A\cdot)$,
$$B_k(y_0) \subset \{v < 1\} \subset B_{nk}(y_0).$$
That is, $A^{-1}E = kB_1$, where $E$ is the John ellipsoid of $\{u < 1\}$.
We claim that $1/n \leq k \leq 2$. Indeed, if the left inequality is false, then $|y-y_0|^2$ must lie strictly below $v$ in $\{v < 1\}$ by the maximum principle, but this is false at the origin. Similarly, if the right inequality is false, then $|y-y_0|^2/4$ must lie strictly above $v$ in $B_2(y_0)$, but this is false at $y_0$. See Figure \ref{MAMaxPr}.

We conclude immediately that $\{v < 1\} \subset B_{4n}(0)$. In addition, the Pogorelov estimate implies that $|D^2v| \leq C(n)$ in $\{v < 1/2\}$. Elementary calculus thus implies that $v \leq C(n)|y|^2$ in $\{v < 1/2\}$, whence $B_{c(n)}(0) \subset \{v < 1\}$. We conclude that 
$$AB_{c(n)} \subset \{u < 1\} \subset AB_{4n}.$$ 
Finally, $|D^2v(0)| = |A^TA| \leq C(n)$ (again by Pogorelov). Using that $\det A = 1$ we conclude that its principal values are bounded between positive dimensional constants, and the desired containments (\ref{MAContainment}) follow.

\begin{rem}
One may ask what happens when $u$ is not convex. In two dimensions the only alternative is that $u$ is concave, so up to taking $-u$ convexity is automatic. In higher dimensions there are a plethora of non-convex global solutions; consider for example $u(x_1,\,x_2,\,x_3) = f(x_1) + x_1x_2 -x_3^2/2$, where $f$ is arbitrary. Two major issues are that convexity is already a quite rigid condition, and losing it opens up many possibilities (note e.g. that convex solutions to $\Delta u = 1$ are quadratic polynomials by the Liouville theorem, while there are many more non-convex solutions), and that the equation is no longer elliptic when $u$ is not convex.
\end{rem}

\begin{rem}
It is worth remarking that rigidity is false for the complex Monge-Amp\`{e}re equation, even in (complex) dimension $n = 2$. Indeed, the function $u(z_1,\,z_2) = 2|z_1|(1+|z_2|^2)$ solves $\det \partial \bar{\partial} u = 1$ on $\mathbb{C}^2$ (\cite{Bl}). One issue is that solutions are not necessarily convex. Another is the lack of a (known) analogue of Pogorelov's interior $C^2$ estimate. Thus, even if one assumes convexity, it is not clear that rigidity should hold. This is related to the fact that the equation is invariant under adding the real part of a holomorphic function, which can destroy convexity.
\end{rem}

\newpage
\section{Minimal Surface Equation in Higher Dimensions}\label{MSEHigherD}
In two dimensions, we could prove the Bernstein theorem for minimal graphs using the equation (first variation of the area) directly. Extending the Bernstein theorem for minimal graphs to higher dimensions requires two key tools. The first is the stability inequality (coming from the second variation of area), which reflects that minimal graphs minimize area. The second is the monotonicity formula, which allows one to reduce the problem to studying cones.

\subsection{Small Variations: Mean Curvature and Stability}
Let $\Sigma$ be a smooth, oriented hypersurface in $\mathbb{R}^{n+1}$ with unit normal $\nu$, second fundamental form $II,$ and mean curvature $H$. Below, $\nabla_{\Sigma}$ and $\Delta_{\Sigma}$ will denote the usual gradient and Laplace operators on $\Sigma$. Given a smooth function $\psi$ on $\Sigma$, these can be calculated at $x \in \Sigma$ by extending $\psi$ near $x$ to be constant in the direction $\nu(x)$, taking the usual gradient and Laplace of the extension in $\mathbb{R}^{n+1}$, and evaluating at $x$.

\begin{figure}
 \centering
    \includegraphics[scale=0.7, trim={0mm 115mm 0mm 30mm}, clip]{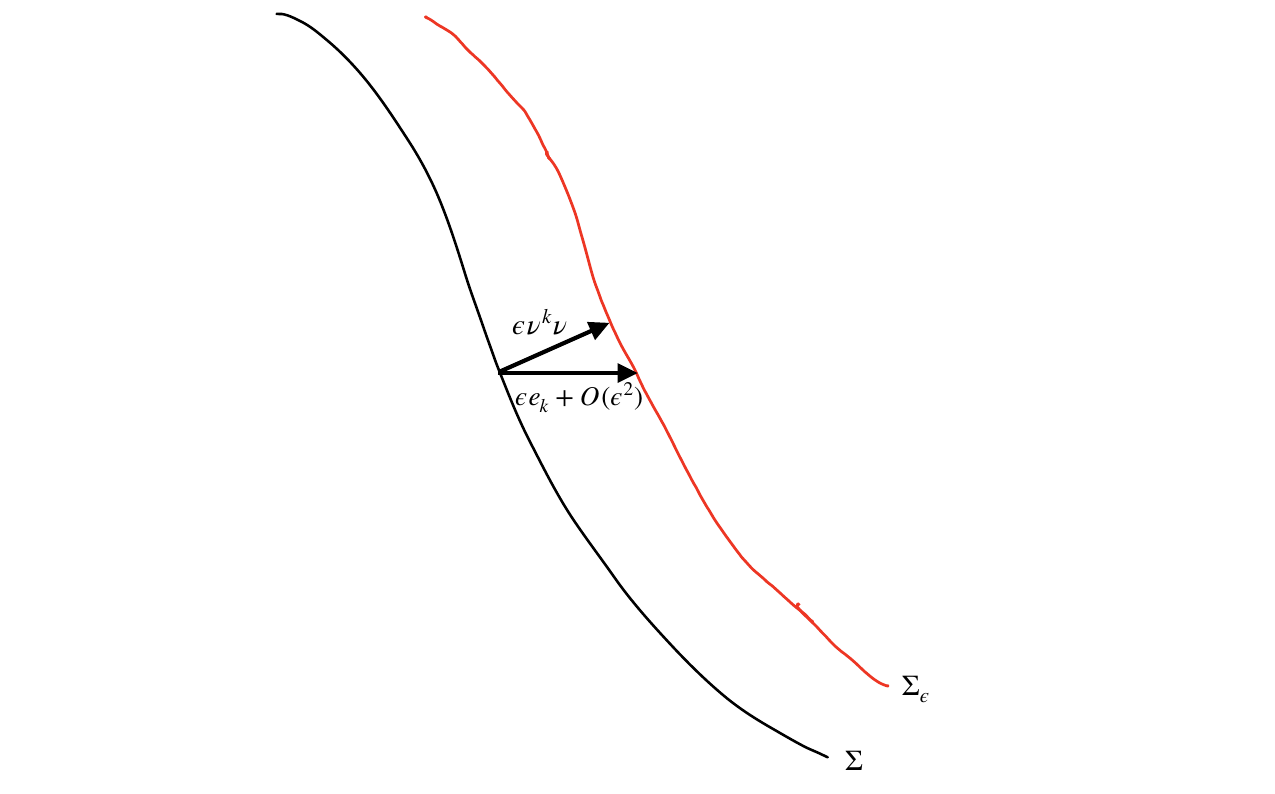}
\caption{}
\label{Transl}
\end{figure}

Let $\varphi$ be a smooth function on $\Sigma$, and for $x \in \Sigma$ and $\epsilon$ small consider the map $T(x) = x + \epsilon \varphi \nu(x)$, giving rise to a perturbed surface $\Sigma_{\epsilon} = T(\Sigma)$ with mean curvature $H_{\epsilon}$. We have the following expansion:
\begin{equation}\label{HChange}
H_{\epsilon}(T(x)) = H(x) + \epsilon(\Delta_{\Sigma} \varphi + c^2\varphi)(x) + O(\epsilon^2).
\end{equation}
Here and below, $c^2$ denotes $|II|^2$, the sum of squares of principal curvatures. We let
$$L := \Delta_{\Sigma} + c^2,$$ 
denote the Jacobi operator. A convenient way to derive (\ref{HChange}) is to write $\Sigma$ locally as a graph over its tangent plane. After a translation and a rotation we may assume that the point of tangency is the origin, and $\nu(0) = e_{n+1}$. Then $\Sigma = \{(y,\,u(y))\}$ locally, with $y \in \mathbb{R}^n$ and $u(0) = |\nabla u(0)| = 0$. Using that 
$$\nu = (-\nabla u,\,1)(1+|\nabla u|^2)^{-1/2} = (-\nabla u,\, 1-|\nabla u|^2/2) + O(|y|^3),$$ 
we see that the perturbed surface $\Sigma_{\epsilon}$ is locally the graph of a function $w$ such that
$$w(y-\epsilon\varphi\nabla u(y)) = u(y) + \epsilon\varphi(1-|\nabla u|^2/2) + O(|y|^3).$$
Here we have extended $\varphi$ to be constant in the $e_{n+1}$ direction. Differentiating this identity (in $y \in \mathbb{R}^n$) and evaluating at the origin twice gives
$$\nabla w(0) = \epsilon \nabla \varphi(0) + O(\epsilon^2), \quad D^2w(0) = D^2u(0) + \epsilon(D^2\varphi(0) + \varphi [D^2u(0)]^2) + O(\epsilon^2).$$
Identity (\ref{HChange}) follows, using that $II(0) = D^2u(0)$ and that $H_{\epsilon}(\epsilon\varphi e_{n+1}) = \Delta w(0) + O(|\nabla w|^2)$. 
As a consequence of (\ref{HChange}), if $\Sigma$ has constant mean curvature and $\Sigma_{\epsilon}$ is an isometry of $\Sigma$ up to an error of size $\epsilon^2$ for all $\epsilon$, then $L\varphi = 0$. Since $\varphi = \nu \cdot e_k := \nu^k$ generates a translation in the direction $e_k$ (see Figure \ref{Transl}), we conclude that if $\Sigma$ has constant mean curvature, then $L\nu^k = 0$ for all $k = 1,\,...,\,n+1$.

We now consider how area changes. We drop the subscript $\Sigma$ from the operators for simplicity. In coordinates where the unit normal is the last coordinate direction, the tangential differential of $T$ is $I - \epsilon \varphi II$ in the top $n$ rows, and $\epsilon \nabla \varphi$ in the bottom row. Thus, the square of the area element of $\Sigma_{\epsilon}$ is
$$\det(DT^TDT) = \det(I - 2\epsilon \varphi II + \epsilon^2(\varphi^2 II^2 + \nabla \varphi \otimes \nabla \varphi)).$$
Expanding the RHS and using that $2\sigma_2(II) = H^2-c^2$, we get that the square of the area element is
$$1 - 2\epsilon H\varphi + \epsilon^2 (|\nabla \varphi|^2 + (2H^2-c^2) \varphi^2) + O(\epsilon^3).$$
We conclude after taking the square root that
$$\text{Area}(\Sigma_{\epsilon}) = \text{Area}(\Sigma) - \epsilon \int_{\Sigma} H\varphi + \frac{\epsilon^2}{2}\int_{\Sigma} (|\nabla \varphi|^2 + (H^2 - c^2)\varphi^2) + O(\epsilon^3).$$
We now restrict to the case that $\Sigma$ is minimal, so that $H = 0$, and we let $\Omega$ be a domain in $\Sigma$. We say that $\Omega$ is stable if small normal variations supported in $\Omega$ cannot decrease area, to second order. That is,
\begin{equation}\label{Stab1}
\int_{\Sigma} (|\nabla \varphi|^2 - c^2 \varphi^2) \geq 0 \text{ for all } \varphi \in C^{\infty}_0(\Omega).
\end{equation}
After integrating by parts this means
\begin{equation}\label{Stab2}
\int_{\Sigma} \varphi L(\varphi) \leq 0 \text{ for all } \varphi \in C^{\infty}_0(\Omega).
\end{equation}
Equivalently, the maximum Dirichlet eigenvalue of $L$ in $\Omega$ is nonpositive.

The latter expression gives some useful criteria for stability and instability. First, if there exists a function $f \in H^1_0(\Omega)$ such that $f \geq 0,\, Lf \geq 0$ in $\Omega$ and $fLf > 0$ somewhere, we immediately get instability. Alternatively, if there is a positive super-solution on some domain in $\Sigma$ containing $\overline{\Omega}$ (that is, $Lf \leq 0$ and $f > 0$), we get stability of $\Omega$. Indeed, if the maximal Dirichlet eigenvalue of $L$ is positive, then corresponding eigenfunction is a positive sub-solution of the Jacobi equation, and we can rule this out using the maximum principle (a multiple of $f$ touches the eigenfunction from above at some point, a contradiction). A consequence is that minimal graphs $\{(x,\,u(x))\}$ are stable, because the vertical component $\nu^{n+1} = (1+|\nabla u|^2)^{-1/2}$ of the upper unit normal is positive and solves $L\nu^{n+1} = 0$.

\subsection{Minimizing Property and Monotonicity Formula: Reduction to Cones}
We now reduce the Bernstein problem to the study of minimal cones. The ideas described in this sub-section are due to Fleming \cite{Flem} and De Giorgi \cite{DG2}, and a detailed development can be found in the book of Giusti \cite{Giu}. In this sub-section $\Sigma$ denotes an entire minimal graph in $\mathbb{R}^{n+1}$, graphical over $\{x_{n+1} = 0\}$.

\begin{figure}
 \centering
    \includegraphics[scale=0.7, trim={0mm 170mm 0mm 30mm}, clip]{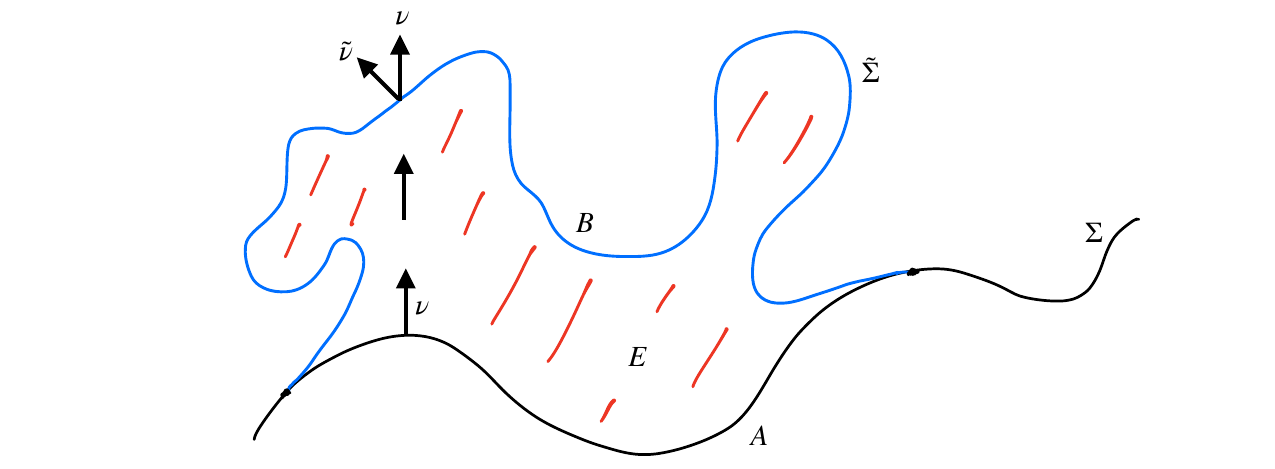}
\caption{}
\label{Calibration}
\end{figure}

The first key observation is that $\Sigma$ is not only stable (see previous subsection), but area minimizing. To see this, extend the upper unit normal $\nu$ vertically in the $x_{n+1}$ direction. The minimal surface equation says that $\text{div}(\nu) = 0$. Let $\tilde{\Sigma}$ be a competitor that agrees with $\Sigma$ outside of a large ball, with unit normal $\tilde{\nu}$, and let $E$ be the region in between $\Sigma$ and $\tilde{\Sigma}$. Then $\partial E$ consists of two portions, $A \subset \Sigma$ and $B \subset \tilde{\Sigma}$ (see Figure \ref{Calibration}). The the divergence theorem says
$$0 = \int_E \text{div}(\nu) = \int_{B} \nu \cdot \tilde{\nu} - \text{Area}(A) \leq \text{Area}(B) - \text{Area}(A)$$
as desired. This is called a calibration argument.

\begin{rem}
There is an efficient proof of the Bernstein theorem in two variables based on the stability inequality and the area-minimizing property of graphs. Take the log cutoff as the normal variation in the stability inequality (\ref{Stab1}). One gets that $c^2 \equiv 0$ provided it took no Dirichlet energy to cut off. This is guaranteed if the graph has quadratic area growth, which holds by the area minimizing property.
\end{rem}

\begin{figure}
 \centering
    \includegraphics[scale=0.7, trim={0mm 135mm 0mm 30mm}, clip]{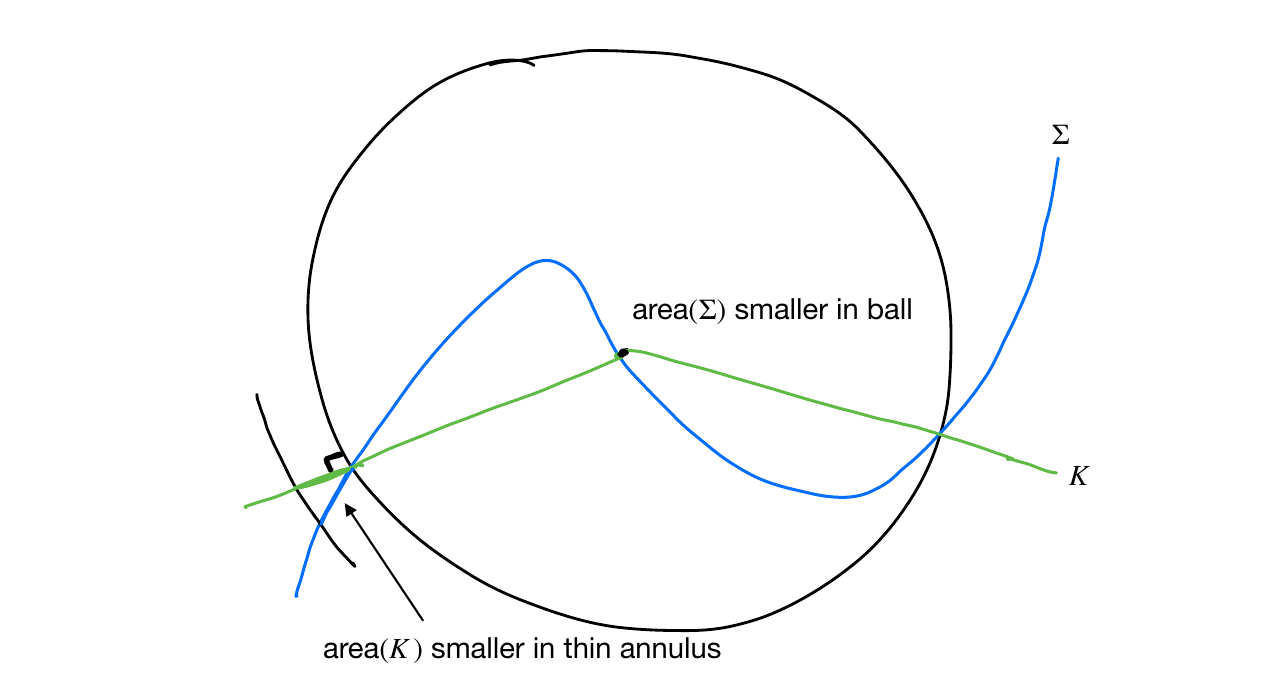}
\caption{}
\label{MF}
\end{figure}

The second key tool is the monotonicity formula. Assume after a translation that $0 \in \Sigma$, and let $B_R$ denote a ball in $\mathbb{R}^{n+1}$ of radius $R$ centered at $0$. A consequence of area minimality is that $V_{\Sigma}(R) := \text{Area}(\Sigma \cap B_R) \leq C(n)R^n$, as can be seen by taking $\partial B_R$ as a competitor. Thus, the quantity 
$$\Phi_{\Sigma}(R) := R^{-n}V_{\Sigma}(R)$$
is bounded. The monotonicity formula says that $\Phi_{\Sigma}$ is non-decreasing in $R$, and constant if and only if $\Sigma$ is a cone, i.e. it is dilation-invariant. The argument goes as follows. We have
$$\Phi_{\Sigma}'(R) = R^{-n-1}(RV'_{\Sigma}(R) - nV_{\Sigma}(R)).$$
Let $K$ be the cone in $\mathbb{R}^{n+1}$ such that $K \cap \partial B_R = \Sigma \cap \partial B_R$. Note that $\Phi_K$ is constant. Euclidean geometry implies that $V_K'(R) \leq V_{\Sigma}'(R)$, with equality only if $\Sigma$ crosses $\partial B_R$ orthogonally. By area minimality, $V_{\Sigma}(R) \leq V_K(R)$. Hence $\Phi_{\Sigma}'(R) \geq V_K'(R) = 0$, and $\Phi_{\Sigma}' \equiv 0$ if and only if $\Sigma$ crosses $B_R$ orthogonally for all $R$, that is, $\Sigma$ is a cone. See Figure \ref{MF}.

Using the previous observations (area-minimizing and monotonicity formula), one can show that for some sequence $R_j \rightarrow \infty$ the blow-down sequence $R_j^{-1}\Sigma$ converges to an area-minimizing cone $K$. The argument goes roughly as follows. The compactness properties of area-minimizing hypersurfaces allow the extraction of an area-minimizing limit hypersurface $K$. Using the monotonicity formula, one can furthermore show that $\Phi_K \equiv \Phi_{\Sigma}(+\infty)$, hence $K$ is a cone. If $K$ can be shown to be a hyperplane, then $\Phi_{\Sigma}(0^+) = \Phi_K = \Phi_{\Sigma}(+\infty)$ (the first equality coming from the fact that $\Sigma$ is locally well-approximated by hyperplanes), hence $\Phi_{\Sigma}$ is constant. Monotonicity implies that $\Sigma$ was the hyperplane $K$ to begin with. In particular, if one proves that area-minimizing cones in $\mathbb{R}^{n+1}$ are hyperplanes, one concludes that entire minimal graphs in $\mathbb{R}^{n+1}$ are hyperplanes.

One can gain an extra dimension by using the graphicality of $\Sigma$ more carefully. It is in fact true that if $\Sigma$ were not a hyperplane, then the cone $K$ is the graph of a function that can take the values $\pm \infty$, and in fact takes the value $+ \infty$ or $-\infty$ on a set in $\mathbb{R}^n$ whose boundary is a non-flat area-minimizing hypercone. Thus, if one proves that area-minimizing hypercones in $\mathbb{R}^n$ are hyperplanes, one concludes that entire minimal graphs in $\mathbb{R}^{n+1}$ are hyperplanes.

To conclude we note that if $K$ is a non-flat area-minimizing hypercone, and it is not smooth away from the origin, then using a blow-up procedure at a point away from the origin one can show that there is a non-flat area-minimizing hypercone in one dimension lower. This is know as Federer's dimension reduction, and it allows one to reduce to the case that $K$ is smooth away from the origin.

\subsection{Simons Inequality}
The formula for the second variation of area indicates that normal variations in regions of high curvature tend to decrease area. It is thus natural to see what happens when we use functions of $c^2$ to test stability. To that end, we calculate the Jacobi operator applied to $c^2$. 

To simplify the calculation, write a smooth minimal hypersurface $\Sigma$ locally as the graph of a function $u$ of $y \in \mathbb{R}^n$ such that $u(0) = |\nabla u(0)| = 0$. The minimal surface equation says that
$$\Delta u = u_{ij}u_iu_j + O(|y|^4),$$
whence
\begin{equation}\label{DerivMSE}
(\nabla \Delta u)(0) = 0, \quad (D^2\Delta u)(0) = 2(D^2u)^3(0).
\end{equation}
We claim that
\begin{equation}\label{cFormula}
c^2 = |D^2u|^2(1-|\nabla u|^2) - 2|D^2u \cdot \nabla u|^2 + O(|y|^4).
\end{equation}
To see this, let $\nu := (-\nabla u,\,1)(1+|\nabla u|^2)^{-1/2}$ be the upper unit normal, let $\nabla$ denote derivatives in $y$, and $\nabla_{\Sigma}$ denote the gradient operator on $\Sigma$. We recall that
$$c^2 = \sum_{k = 1}^{n+1} |\nabla_{\Sigma} \nu^k|^2 = \sum_{k = 1}^{n+1} |\nabla \nu^k - (\nabla \nu^k \cdot \nu)\nu|^2 = \sum_{k = 1}^{n+1} (|\nabla \nu^k|^2 - (\nabla \nu^k \cdot \nu)^2).$$
Using the definition of $\nu$ we calculate:
$$\nu^k = -u_k(1-|\nabla u|^2/2) + O(|y|^5),\, k \leq n, \quad \nu^{n+1} = 1-|\nabla u|^2/2 + O(|y|^4)$$
$$\partial_i\nu^k = -u_{ik}(1-|\nabla u|^2/2) + u_{il}u_lu_k + O(|y|^4),\, k \leq n, \quad \partial_i\nu^{n+1} = -u_{il}u_l + O(|y|^3)$$
$$\nabla \nu^k \cdot \nu = u_{ik}u_i + O(|y|^3),\, k \leq n, \quad \nabla \nu^{n+1} \cdot \nu = O(|y|^2).$$
Plugging these expressions into the formula for $c^2$ yields (\ref{cFormula}).

We now claim that
\begin{equation}\label{CoordSimons}
\frac{1}{2}\Delta c^2(0) = \sum_{i,\,k,\,l} u_{ikl}^2(0) - c^4(0).
\end{equation}
This is the basic Simons identity, written in coordinate-independent fashion as
\begin{equation}\label{SimonsIdentity}
\frac{1}{2}\Delta_{\Sigma} c^2 = |\nabla_{\Sigma} II|^2 - c^4.
\end{equation}
Indeed, directly taking the Laplacian of the expression for $c^2/2$ and evaluating at zero gives the RHS in (\ref{CoordSimons}), plus $\text{tr}[D^2uD^2\Delta u - 2(D^2u)^4].$ Since $D^2uD^2\Delta u = 2(D^2u)^4$ at $0$ by (\ref{DerivMSE}), this vanishes, proving the claim. 

Another way to write (\ref{CoordSimons}) is
\begin{equation}\label{SimonsIdentity2}
cL(c)(0) = \sum_{i,j,k} u_{ijk}(0)^2 - |\nabla c|^2(0).
\end{equation}
We wish to bound the last term on the RHS. By (\ref{cFormula}) we have $c^2 = |D^2u|^2 + O(|y|^2)$. Choosing coordinates in which $D^2u(0)$ is diagonalized and differentiating gives
$$\partial_jc(0) = \left(\sum_{i} u_{ii}u_{iij}/c\right)(0),$$
hence by Cauchy-Schwarz,
$$|\nabla c|^2(0) \leq \sum_{i,\,j} u_{iij}^2(0).$$
Using this in the identity (\ref{SimonsIdentity2}) gives, at $0$,
\begin{equation}\label{SimonsIdentity3}
cL(c) \geq \sum_{i,j,k} u_{ijk}^2 - \sum_{i,j}u_{iij}^2 \geq \left(\sum_i u_{iii}^2 + 3\sum_{i \neq j} u_{iij}^2\right) - \sum_{i,j} u_{iij}^2 = 2\sum_{i \neq j} u_{iij}^2.
\end{equation}
Thus $c$ is a nonnegative sub-solution of the Jacobi equation, which is promising for proving instability results.

To proceed we specialize to the case that $\Sigma$ is a minimal cone. We may assume after a rotation that the vertex is at the origin, and that on the positive $x_n$ axis, $\Sigma$ is locally the graph in the $x_{n+1}$ direction of a one-homogeneous function $u$ whose value and gradient vanish on the positive $x_n$ axis. After rotating in the $x_1,\,...,\,x_{n-1}$ variables we may assume that $D^2u$ is diagonalized along the axis. Then along the $x_n$ axis, (\ref{SimonsIdentity3}) implies
$$cL(c) \geq 2\sum_{i < n} u_{iin}^2.$$
Since $u$ is one-homogeneous, along the axis, we have that $u_{iin} = -u_{ii}/x_n$, hence the RHS can be written $2c^2/|x|^2$. We have arrived at the cone Simons inequality
\begin{equation}\label{ConeSimons}
L(c) \geq 2c/|x|^2,
\end{equation}
again promising for proving instability results.


\subsection{Simons Theorem}
We now prove the celebrated Simons theorem \cite{Si} that stable minimal hypercones in $\mathbb{R}^{n+1}$ are flat when $n \leq 6$, implying the Bernstein theorem (entire minimal graphs in $\mathbb{R}^{n+1}$ are flat) up to dimension $n = 7$. 

\begin{figure}
 \centering
    \includegraphics[scale=0.7, trim={0mm 180mm 0mm 30mm}, clip]{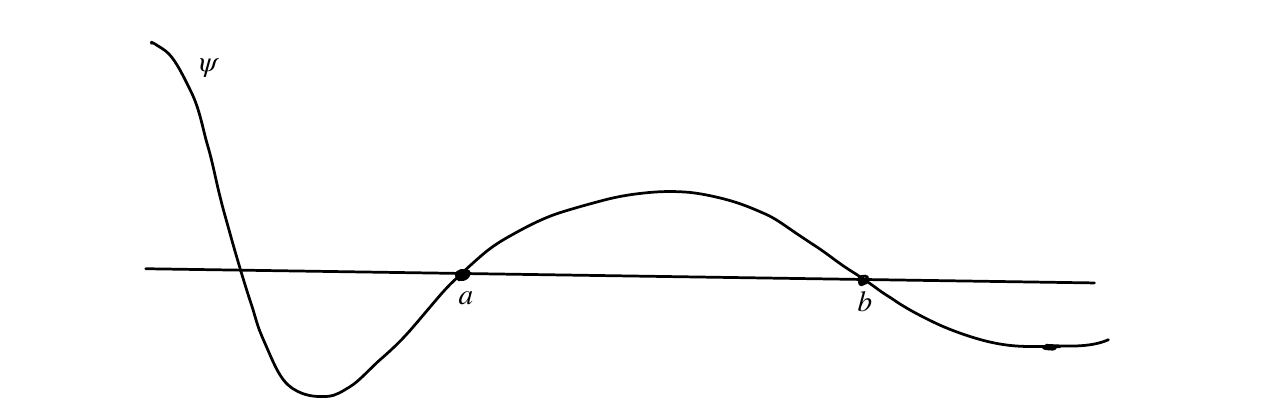}
\caption{}
\label{Oscillate}
\end{figure}

Let $\Sigma$ be a stable minimal hypercone in $\mathbb{R}^{n+1}$, smooth outside of the origin. Let $r := |x|$ and let $\psi$ be a radial function. We have
$$L(c\psi) = L(c)\psi + 2\psi'\partial_rc + c\Delta_{\Sigma}\psi.$$
Homogeneity gives $\partial_rc = -c/r$, and for radial $\psi$ we have $\Delta_{\Sigma}\psi = \psi'' + (n-1)r^{-1}\psi'$. Using these identities along with (\ref{ConeSimons}) we arrive at
$$L(c\psi) \geq c(\psi'' + (n-3)r^{-1}\psi' + 2r^{-2}\psi).$$
The ODE $g'' + (n-3)r^{-1}g' + \gamma r^{-2}g = 0$ has solutions of the form $r^{\lambda}$, where
$$\lambda^2 + (n-4)\lambda + \gamma = 0.$$
We can guarantee that a solution has two zeros $0 < a < b$ and is positive for $r \in (a,\,b)$ provided
$(n-4)^2 < 4\gamma$, that is, $n < 4 + 2\sqrt{\gamma}$ (see Figure \ref{Oscillate}). If we fix $\gamma \in (1,\,2)$, this is satisfied when $n \leq 6$. Taking $\psi$ to be said solution, we get 
$$L(c\psi) \geq (2-\gamma)r^{-2}c\psi.$$
Thus, on the domain $\Omega$ given by the intersection of $\Sigma$ with the annulus $B_b\backslash B_a$, the function $f = c\psi$ is a nonnegative sub-solution to the Jacobi equation that vanishes on $\partial \Omega$, and at any point where $c > 0$, $fL(f) > 0$. As a consequence, if $\Sigma$ is not a hyperplane and $n \leq 6$, then $\Sigma$ cannot be stable, giving the Simons theorem.

\begin{rem}
Almgren \cite{Alm} showed that stable minimal hypercones in $\mathbb{R}^4$ are flat by proving two facts about minimal hypercones $K \subset \mathbb{R}^4$: (1) If the link $K \cap \mathbb{S}^3$ doesn't have the topology of $\mathbb{S}^2$, then $K$ is unstable, and (2) if $K \cap \mathbb{S}^3$ has the topology of $\mathbb{S}^2$, then $K$ is flat.
\end{rem}

\begin{rem}\label{ConditionalBernstein}
One can extend the Bernstein theorem to all dimensions with growth hypotheses. In Section \ref{UnifElliptic} we saw that bounded gradient suffices. In \cite{EH}, Ecker and Huisken relaxed this to sub-linear gradient growth. This is essentially optimal, as in high dimension $n$, the Bombieri-De Giorgi-Giusti entire minimal graphs (see Section \ref{BDGEx}) have gradient growing at the rate $|x|^{1 + O(1/n)}$.
\end{rem}

\subsection{Minimal Surface System}\label{MSSBernstein}
As seen above, there are many entire solutions to the minimal surface system, even when $n = m = 2$ (e.g. any holomorphic map). It is natural to ask for rigidity theorems under the additional hypothesis that the gradient is bounded, since this guarantees linearity in codimension one.

The form of the system (\ref{MSS}) and the discussion in Remark \ref{2D} show that Lipschitz entire solutions are linear for $n = 2$, $m$ arbitrary. To go to higher dimension we require the monotonicity formula, which we proved above for area-minimizing hypersurfaces, but in fact holds in greater generality: For any smooth minimal (not necessarily minimizing) submanifold $\Sigma$ of dimension $n$ in $\mathbb{R}^{n + m}$ containing the origin, the quantity $R^{-n}\text{Area}(\Sigma \cap B_R)$ is non-decreasing in $R$, and it is constant if and only if $\Sigma$ is a cone. Here $B_R$ is an extrinsic ball of radius $R$ centered at $0$. Using the monotonicity formula and a blow-down argument similar to that outlined above, the linearity of global Lipschitz solutions to (\ref{MSS}) would follow from the linearity of one-homogeneous solutions that are smooth outside of the origin. Lipschitz entire solutions to (\ref{MSS}) are thus linear when $n = 3$, by the discussion in Remark \ref{2D}.

The Lawson-Osserman example (see Section \ref{LOEx} below) shows that when $n \geq 4$, there are non-flat graphical minimal cones, so additional conditions are required in higher dimension. A sufficient condition for the linearity of Lipschitz entire solutions to (\ref{MSS}), discovered by M.-T. Wang \cite{Wa}, is the area-decreasing condition: for some $\delta > 0$ and all $x$, the principal values $\lambda_i$ of $Du(x)$ satisfy that $\lambda_i\lambda_j \leq 1-\delta$ for all $i \neq j$. 
We note that the Lawson-Osserman cone has codimension three, so it is feasible that one has stronger results in codimension two.

Finally, we recall that stability played an important role in the proof of the Bernstein theorem for minimal graphs of codimension one. In higher codimension, minimal graphs are not necessarily stable, and moreover, the role of stability is more mysterious. There are interesting results in the case $n = 2$. For example, it can be shown that complete, oriented, stable minimal surfaces of dimension two in $\mathbb{R}^{2+m}$ are contained in an even-dimensional affine subspace and holomorphic with respect to some complex structure, under some additional assumptions e.g. about area growth, topology, and/or total curvature, see the work of Micallef \cite{Mic} and Fraser-Schoen \cite{FS}.

\newpage
\section{Nonlinear Global Solutions}\label{BDGEx}
In this section we build the Bombieri-De Giorgi-Giusti example of an entire minimal graph in $\mathbb{R}^9$ \cite{BDG}. We follow the approach taken in \cite{MY} to solve the Bernstein problem for anisotropic minimal surfaces.

\subsection{Foliation}
To begin we observe that the Simons cone 
$$S := \{|x| = |y|\} \subset \mathbb{R}^8,\, x,\,y \in \mathbb{R}^4,$$ 
which is minimal away from the origin, is stable. Indeed, using that $c = 6/r^2$ (here $r^2 = |x|^2 + |y|^2$) it is easy to verify that
$$Lr^{-2} = Lr^{-3} = 0,$$
hence there are positive solutions to the Jacobi equation. (To be precise, this shows that $S \backslash B_{r}$ is stable for all $r > 0$. One can show global stability by cutting off variations near the vertex and using that area in $B_{r}$ scales like $r^{7}$).

\begin{figure}
 \centering
    \includegraphics[scale=0.65, trim={0mm 150mm 0mm 15mm}, clip]{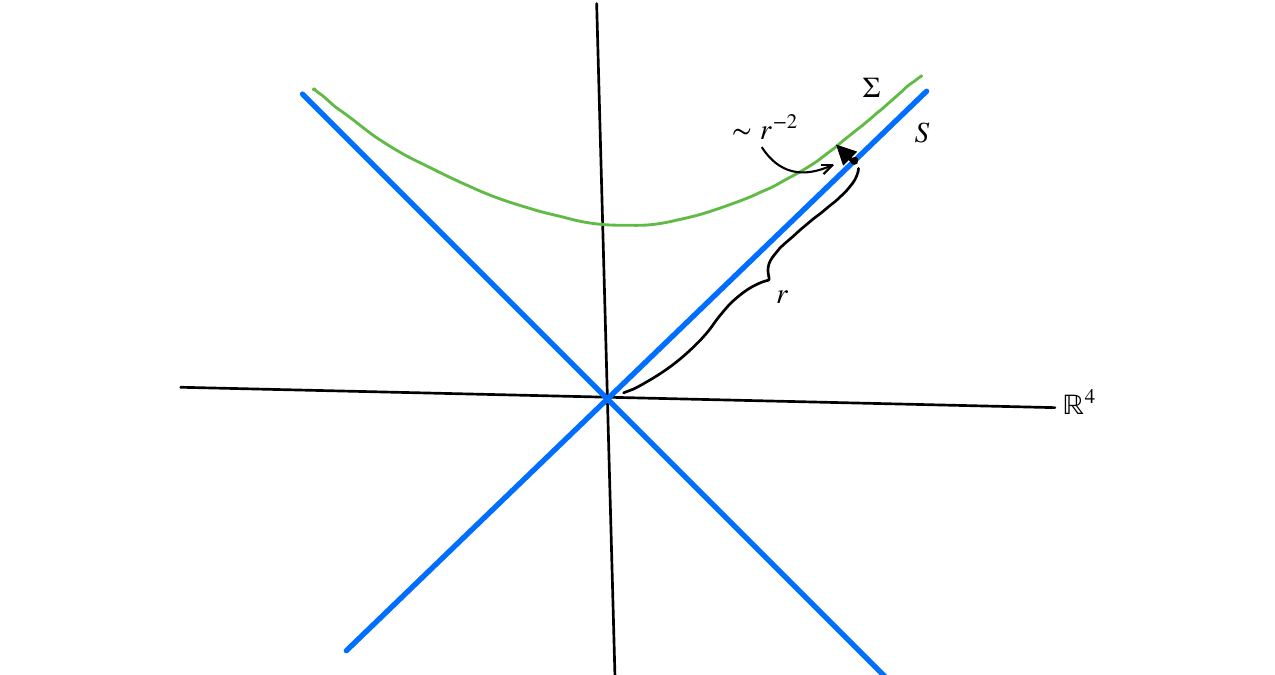}
\caption{}
\label{Leaf}
\end{figure}

These Jacobi fields suggest the existence of minimal surfaces close to $S$ that lie on one side of $S$ (in contrast with lower dimensions, where the Simons theorem says that minimal surfaces nearby non-flat minimal cones oscillate around the cones). We confirm this now. The main claim of this sub-section is that there is a smooth, even, locally uniformly convex function $\sigma$ on $\mathbb{R}$ such that $\sigma(0) = 1$, for $s > 1$ and some $a > 0$ we have
$$\sigma(s) = s + as^{-2} + O_2(s^{-3}),$$
and furthermore
$$\Sigma := \{|y| = \sigma(|x|)\}$$
is minimal (see Figure \ref{Leaf}). Here and below, given a function $f$ on $(1,\,\infty)$, $O_2(f)$ denotes a function whose value, derivative, and second derivative are bounded by a constant times those of $f$ on $(1,\,\infty)$.

\begin{figure}
 \centering
    \includegraphics[scale=0.65, trim={0mm 130mm 0mm 15mm}, clip]{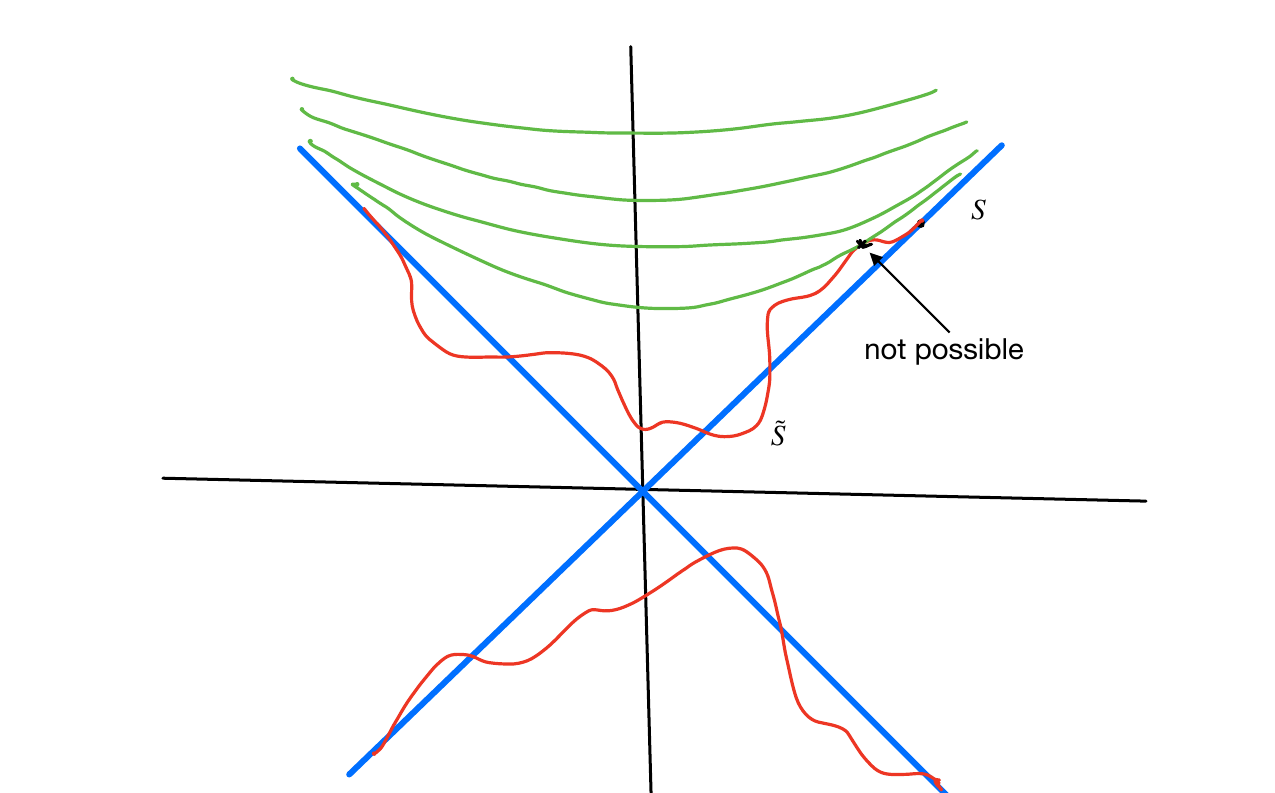}
\caption{}
\label{Minimizing}
\end{figure}

The properties of $\sigma$ imply that the dilations $\lambda \Sigma,\, \lambda > 0$ foliate one side of $S$. By symmetry, each side of $S$ is foliated by smooth minimal hypersurfaces approaching at the same rate as the first Jacobi field ($r^{-2}$). The foliation in fact implies that $S$ is area-minimizing. Roughly, if solving the Plateau problem with the same boundary as $S$ on some domain gave a different hypersurface $\tilde{S}$, then a leaf in the foliation would touch $\tilde{S}$ from one side, a contradiction of the strong maximum principle (see Figure \ref{Minimizing}).

\begin{figure}
 \centering
    \includegraphics[scale=0.7, trim={0mm 160mm 0mm 15mm}, clip]{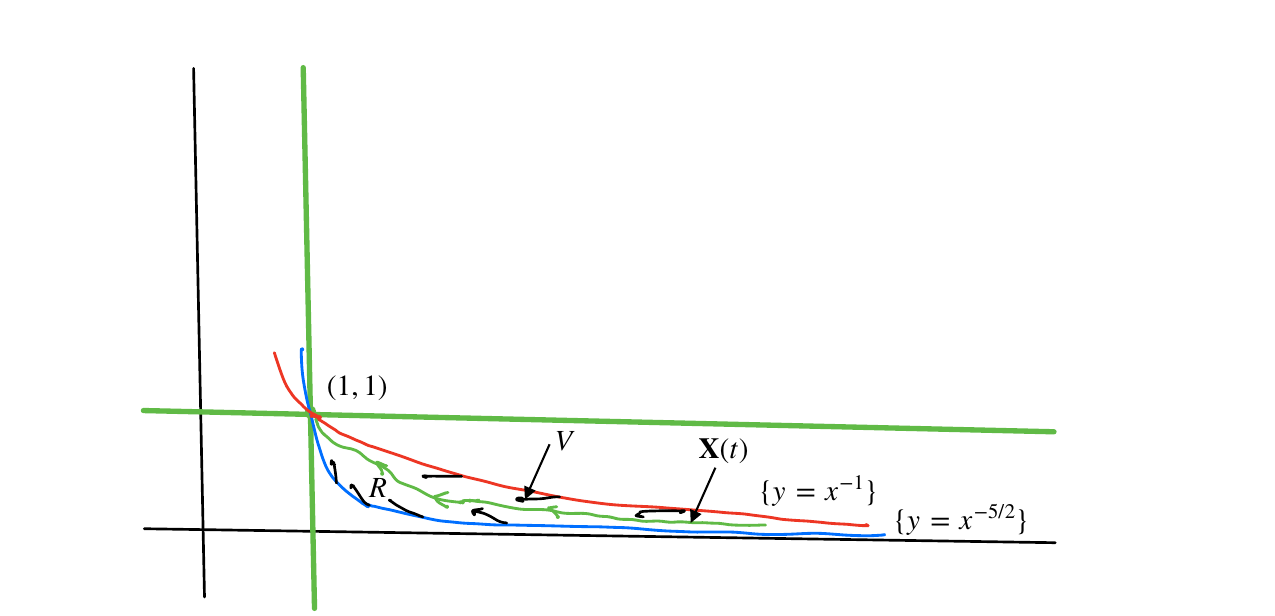}
\caption{}
\label{TrappingFig1}
\end{figure}

To prove the claim we need to solve the ODE
\begin{equation}\label{ODE}
G(\sigma) := \sigma''(s) + 3(1+\sigma'^2(s))\left(\frac{\sigma'(s)}{s} - \frac{1}{\sigma(s)}\right) = 0.
\end{equation}
Indeed, the mean curvature of $\{|y| = \sigma(|x|)\}$ is given by $(1+\sigma'^2)^{-3/2}G(\sigma)$. This can be derived geometrically (the first term corresponds to the curvature coming from bending in the graph of $\sigma$, and the second two correspond to the curvatures coming from rotations around the vertical, resp. horizontal copies of $\mathbb{R}^4$), or by taking the first variation of the area $\text{const.}\int \sqrt{1+\sigma'^2}\sigma^3s^3\,ds$. The local solvability  of (\ref{ODE}) near $0$ for a uniformly convex, even solution with $\sigma(0) = 1$ is standard.
Letting ${\bf X}(t) = (e^{-t}\sigma(e^t),\,\sigma'(e^t))$, we see that (\ref{ODE}) is equivalent to the autonomous system
$${\bf X}'(t) = V({\bf X}),$$
where
$$V(x,\,y) = \left(-x + y,\,3(1+y^2)\left(\frac{1}{x} - y\right)\right).$$
We note that $(1,\,1)$ is a zero of $V$. We let
$$R := \{x \geq 1\} \cap \{x^{-5/2} \leq y \leq x^{-1}\}.$$
It is not hard to calculate that on $\partial R \backslash \{(1,\,1)\}$, the vector field $V$ points into $R$ (see Figure \ref{TrappingFig1}). Indeed, along the top curve, the vertical component of $V$ vanishes, while the horizontal one is negative. For the bottom curve, one needs to check that
$$\frac{3(1+x^{-5})(x^{-1}- x^{-5/2})}{x-x^{-5/2}} > \frac{5}{2}x^{-7/2},\, x > 1,$$
which reduces after some manipulations to proving that
$$P(z) = 6z^{13} - 11z^{10} + 11z^{3} - 6 > 0,\, z > 1.$$
This in turn follows from $P(1) = 0,\, P'(1) = 1,\, (P'(z)/z^2)' = 780z^9-770z^6$. Moreover, for $t$ very negative, ${\bf X}$ lies in $R$. Indeed, the uniform convexity of $\sigma$ implies that the second component of ${\bf X}$ is increasing, implying that ${\bf X}$ lies below the top curve of $\partial R$. It also implies that $\sigma\sigma'(s)/s \geq c > 0$ for $s$ small, hence the solution curve lies above $y = cx^{-1} > x^{-5/2}$ for $x$ large. Since $V$ has negative first component in $R$ we conclude that ${\bf X}$ is trapped in $R$ and tends to $(1,\,1)$ as $t$ tends to infinity.

\begin{figure}
 \centering
    \includegraphics[scale=0.7, trim={0mm 155mm 0mm 15mm}, clip]{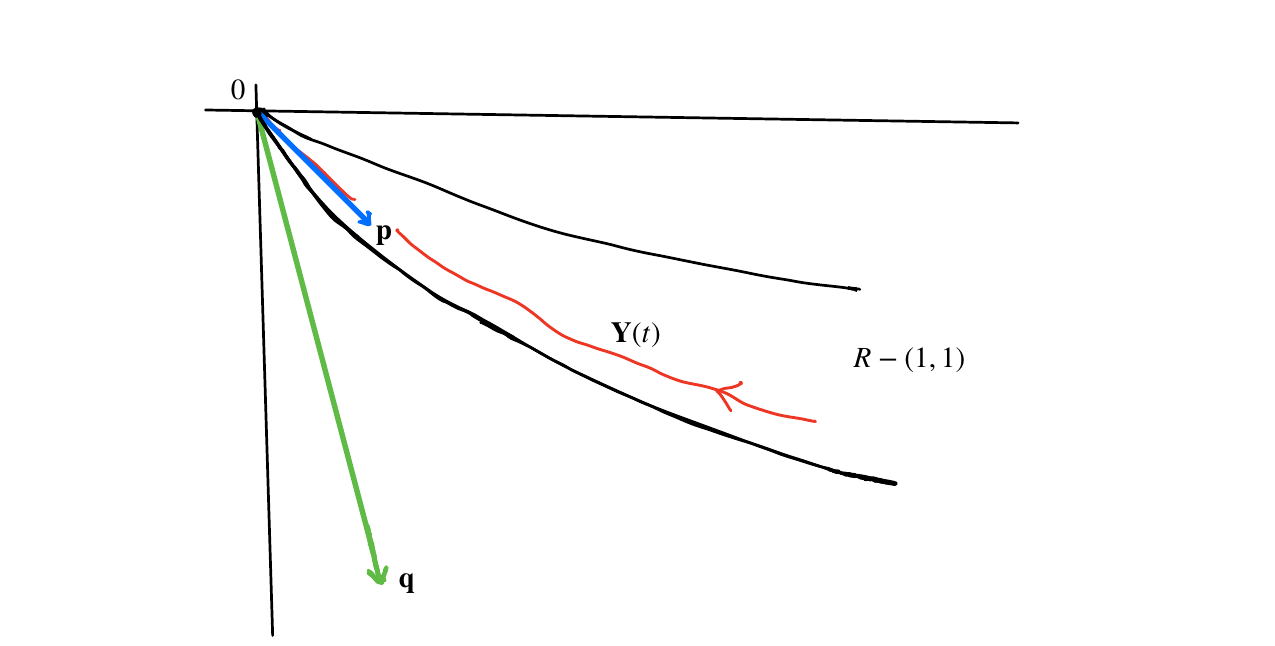}
\caption{}
\label{TrappingFig2}
\end{figure}

The rest is an analysis of the linearized problem at $(1,\,1)$. Let ${\bf X} = {\bf Y} + (1,\,1)$, so that ${\bf Y}$ tends to zero, and is trapped between (translated) boundary curves of $R$, which we note have slopes $-5/2,\,-1$ at $0$. Expanding $V$ at $(1,\,1)$ gives 
\begin{equation}\label{YODE}
{\bf Y}' = M{\bf Y} + O(|{\bf Y}|^2), \quad M = \left(\begin{array}{cc}
-1 & 1 \\
-6 & -6
\end{array}\right).
\end{equation}
The matrix $M$ has eigenvectors ${\bf p} := (1,\,-2)$ and ${\bf q} := (1,\,-3)$ corresponding to eigenvalues $-3,\,-4$ respectively. The second components of ${\bf p}$ and ${\bf q}$ reflect the decay rates of the Jacobi fields mentioned above. We assert that, for some $a > 0$ and $b \in \mathbb{R}$,
\begin{equation}\label{YExp}
{\bf Y}(t) = ae^{-3t}{\bf p} + be^{-4t}{\bf q} + O(e^{-6t}).
\end{equation}
Rewriting this in terms of $\sigma$ and using the ODE gives the claim. To show (\ref{YExp}) we first note that in the directions $e$ of lines with slope in $(-5/2,\,-1)$ we have
$$2Me \cdot e \in (-15/2,\,-2).$$
Indeed, for vectors of the form $e = (1,\,-s)/\sqrt{1+s^2}$ we have
$$-2Me \cdot e =  2\frac{6s^2-5s+1}{1+s^2}.$$
The quantity on the right is increasing for $s \geq 1$, and its values at $s = 1$ and $s = 5/2$ are $2,\,\frac{208}{29} \left(< \frac{15}{2}\right)$.
It follows easily using (\ref{YODE}) and the geometry of the trapping region $R$ that, for some $c > 0$ and all $t > 0$,
$$ce^{-4t} \leq |{\bf Y}| \leq c^{-1} e^{-t/2}.$$
To conclude, write ${\bf Y} = a(t){\bf p} + b(t){\bf q}$. Using (\ref{YODE}) and this decomposition one can boost the decay rate to $|{\bf Y}| \leq c^{-1}e^{-3t}$ for some $c > 0$ and all $t > 0$, and then using (\ref{YODE}) once more gives (\ref{YExp}) for some $a,\,b \in \mathbb{R}$. The curves that trap ${\bf Y}$ force $a \geq 0$, and moreover give that $b = 0$ if $a$ is (note that ${\bf q}$ is outside of the trapping region, see Figure \ref{TrappingFig2}). In the latter case the decay rate of $|{\bf Y}|$ is $e^{-6t}$, violating the lower bound of $e^{-4t}$ determined above and completing the construction.

\subsection{Approach to Entire Graphs}
To go from minimal foliation in $\mathbb{R}^8$ to entire minimal graph in $\mathbb{R}^9$, the idea is to pick a function on $\mathbb{R}^8$ whose level sets are the leaves in the foliation. The symmetries of the foliation suggest taking a function that depends only on $|x|$ and $|y|$, vanishes on $S$, and is odd under exchange of $|x|$ and $|y|$. Moreover, the $r^{-2}$ approach rate of the leaves to $S$ at infinity implies at least quadratic gradient growth of such a choice, so it is natural to pick a function that is homogeneous of degree $3$.

Unsurprisingly, picking the $3$-homogeneous function with the above symmetries and $\Sigma$ as its $1$-level set doesn't quite work. Indeed, minimality of the level sets implies that the mean curvature of the graph is a positive multiple of two derivatives of the function in the direction of the gradient (i.e. the infinity Laplacian, see Subsection \ref{SubSuper} for a more precise expression), which does not vanish. However, by perturbing the leaves to have mean curvature of a desired sign and the same symmetries and asymptotics as before, it turns out we can build sub- and super-solutions with cubic growth. We can then use these to ``trap" the exact solution.

\begin{figure}
 \centering
    \includegraphics[scale=0.65, trim={0mm 155mm 0mm 10mm}, clip]{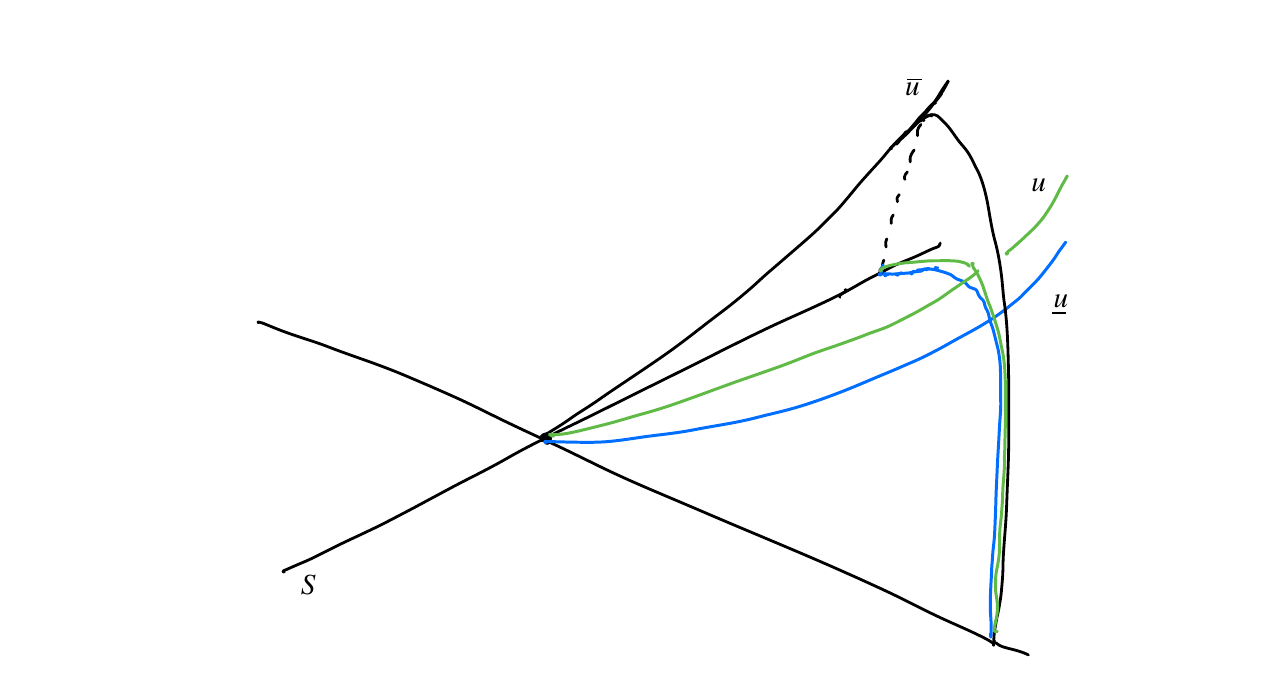}
\caption{}
\label{BDGTrapping}
\end{figure}

More precisely, one builds functions $\underline{u},\,\overline{u}$ with the above symmetries that have the perturbed leaves as level sets, grow cubically at infinity, satisfy $0 \leq \underline{u} \leq \overline{u}$ in $\{|y| > |x|\}$, such that the mean curvature vectors of the graphs of $\underline{u},$ resp. $\overline{u}$ have positive, resp. negative vertical component over $\{|y| > |x|\}$. The solutions to the Dirichlet problem for the MSE in $B_R$ with boundary data $\underline{u}$ then vanish on $S$ by the symmetries of the boundary data, thus lie between $\underline{u},\,\overline{u}$ on each side of $S$ by the maximum principle. By taking $R$ to infinity one obtains the Bombieri-De Giorgi-Giusti graph, which has cubic growth, in the limit (see Figure \ref{BDGTrapping}). 

\begin{rem}
Proving the convergence of the solutions to the Dirichlet problem involves the use of some deep results, notably the interior gradient estimate of Bombieri-De Giorgi-Miranda \cite{BDM}.
\end{rem}

In the following two sub-sections we will show how to construct $\overline{u}$; constructing $\underline{u}$ is very similar (see Remark \ref{SubSoln}).

\subsection{Perturbed Leaves}
\begin{figure}
 \centering
    \includegraphics[scale=0.7, trim={0mm 155mm 0mm 10mm}, clip]{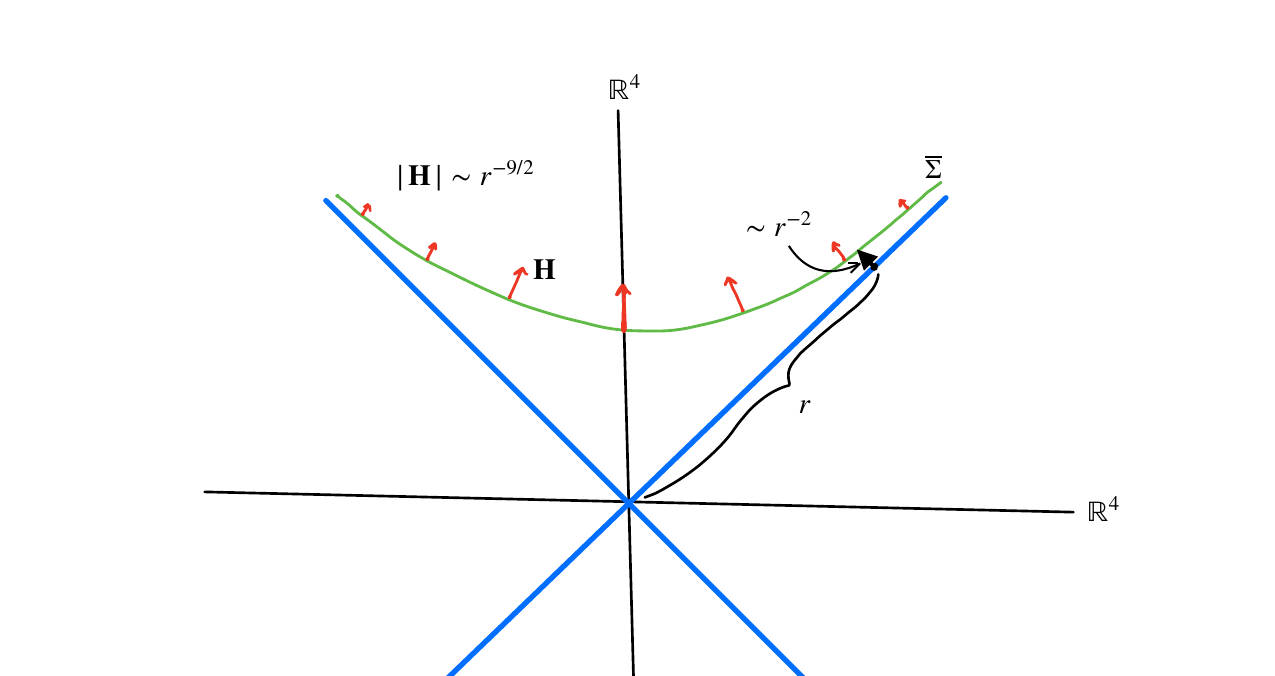}
\caption{}
\label{PertLeaf}
\end{figure}

We now aim to perturb $\Sigma$ to get nearby surfaces with the same asymptotics but mean curvature of a desired sign. We claim that there exists $\overline{\sigma}$ even, smooth, locally uniformly convex, such that for some $\overline{a},\,c > 0$,
$$\overline{\sigma} = s + \overline{a}s^{-2} + O_2(s^{-5/2}),\, s > 1,\, \text{ and }$$
$$G(\overline{\sigma}) \geq c\sigma^{-9/2}.$$
That is, the mean curvature vector of $\overline{\Sigma}: = \{|y| = \overline{\sigma}(|x|)\}$ points away from $S$ and has size decreasing like distance from origin to the power $-9/2$ (see Figure \ref{PertLeaf}).

To do this we study linearized operator $\mathcal{L}$ of $G$ at the solution $\sigma$, given by
$$\mathcal{L}f = f'' + (\log p)'f' + qf,$$
where
$$p = s^3\sigma^3(1+\sigma'^2)^{-3/2}, \quad q = 3(1+\sigma'^2)\sigma^{-2}.$$
One solution $f_0$ to $\mathcal{L}f_0 = 0$ comes from the invariance of the equation $G = 0$ under Lipschitz rescalings:
$$f_0 := -\frac{d}{d\lambda}|_{\lambda = 1} (\lambda^{-1}\sigma(\lambda \cdot)) = \sigma - s\sigma' = 3as^{-2} + O_2(s^{-3}).$$
By writing a solution to $\mathcal{L}f = g$ as the product of $f_0$ with another function, it is not hard to derive the formula for a solution:
$$f(s) = f_0(s)\int_{0}^s f_0^{-2}(t)p^{-1}(t)\int_0^{t} f_0(\tau)p(\tau)g(\tau)d\tau\,dt.$$
Taking $g = \sigma^{-9/2}$ and using the asymptotics of $\sigma$ we get a smooth even solution $f$ such that, for some $c \in \mathbb{R}$,
$$f = cs^{-2} + O_2(s^{-5/2}), \quad s > 1.$$
From here, using Taylor expansion it is not hard to show that
$$|G(\sigma + \epsilon f) - \epsilon \mathcal{L}(f)| \leq C\epsilon^2\sigma^{-7}$$
for some constant $C$, hence
$$G(\sigma + \epsilon_0 f) \geq \frac{1}{2}\epsilon_0 \sigma^{-9/2}$$
for $\epsilon_0 > 0$ small. Taking $\overline{\sigma} = \sigma + \epsilon_0 f$ thus does the job.

\subsection{Supersolution}\label{SubSuper}
We let 
$$w(\lambda s,\,\lambda \bar{\sigma}(s)) = \lambda^3$$
for $s \in \mathbb{R}$ and $\lambda > 0$, and extend $w$ by odd reflection over the diagonals to all of $\mathbb{R}^2$. Let
$$v(x,\,y) = w(|x|,\,|y|).$$
That is, $v$ is the $3$-homogeneous function with the desired symmetries and $1$-level set $\overline{\Sigma}$. On $\{|x| = \lambda s,\, |y| = \lambda \bar{\sigma}(s)\}$ one can show using the asymptotics of $\sigma,\,\overline{\sigma}$ that, for some fixed $C > 0$,
\begin{equation}\label{FunctionEstimates}
v = \lambda^3,\, C^{-1}\lambda^2\sigma(s)^2 \leq |\nabla v| \leq C \lambda^2\sigma(s)^2, \, |D^2v| \leq C \lambda \sigma(s)^{7/2}.
\end{equation}

We will choose $\overline{u}$ to have the same level sets as $v$. More precisely, we will take $\overline{u} = F(v)$, where $F$ is odd, increasing, and has linear growth at infinity, so that $\overline{u}$ has cubic growth. Let $\nu := \nabla v / |\nabla v|$ in $\{|y| > |x|\}$. The condition that $\overline{u}$ is a super-solution to the minimal surface equation in $\{|y| > |x|\}$ is
$$\frac{\overline{u}_{\nu\nu}}{1+\overline{u}_{\nu}^2} - H\overline{u}_{\nu} \leq 0,$$
where $H$ is the mean curvature of the level set with respect to the choice of unit normal $\nu$. In terms of $v$ and $F$,
$$v_{\nu\nu} \leq Hv_{\nu}(1+F'^2(v)v_{\nu}^2) - (F''/F')(v) v_{\nu}^2.$$
Evaluating on the level set $\{|x| = \lambda s,\, |y| = \lambda \bar{\sigma}(s)\}$ and using (\ref{FunctionEstimates}) along with the fact that $H \geq c\lambda^{-1}\sigma^{-9/2}$ on this surface by the considerations in the previous sub-section, we see that this condition is satisfied provided
$$C \leq F'^2(\lambda^3)\lambda^4\sigma^{-2} - (F''/F')(\lambda^3)\lambda^3\sigma^{1/2}$$
for some $C > 0$ fixed and all $\lambda > 0,\, \sigma \geq 1$. Writing $s = \lambda^3,\, t = \sigma^{1/2}$, this becomes
\begin{equation}\label{FinalIneq}
C \leq F'^2(s)s^{4/3}t^{-4} - s(F''/F')(s)t \quad \text{ for all } s > 0,\, t \geq 1.
\end{equation}
We claim that the choice of $F$ determined by
$$F(0) = 0,\, F'(s) = A|s|^{-5/6} + e^{A^2\int_{|s|}^{\infty} \frac{\tau^{-11/12}}{1+\tau^{1/6}}\,d\tau}$$
does the job, provided $A$ is sufficiently large. To verify, split into two cases. When $0 < s \leq 1$, use that $F'^2(s) \geq A^2s^{-5/3}$ and $-s(F''/F')(s) \geq s^{1/12}/2$, so the RHS of (\ref{FinalIneq}) is bounded below by $A^2(s^{1/12}t)^{-4} + (s^{1/12}t)/2$. When $1 < s$, use that $F'^2(s) \geq 1$ and that $-s(F''/F')(s) \geq As^{-1/12}/4$, so the RHS of (\ref{FinalIneq}) is bounded below by $s(s^{-1/12}t)^{-4} + A(s^{-1/12}t)/4 \geq (s^{-1/12}t)^{-4} + A(s^{-1/12}t)/4$. In either case, the minimum possible value of the RHS is a positive power of $A$, hence (\ref{FinalIneq}) is satisfied for $A$ large.

\begin{figure}
 \centering
    \includegraphics[scale=0.7, trim={0mm 75mm 0mm 10mm}, clip]{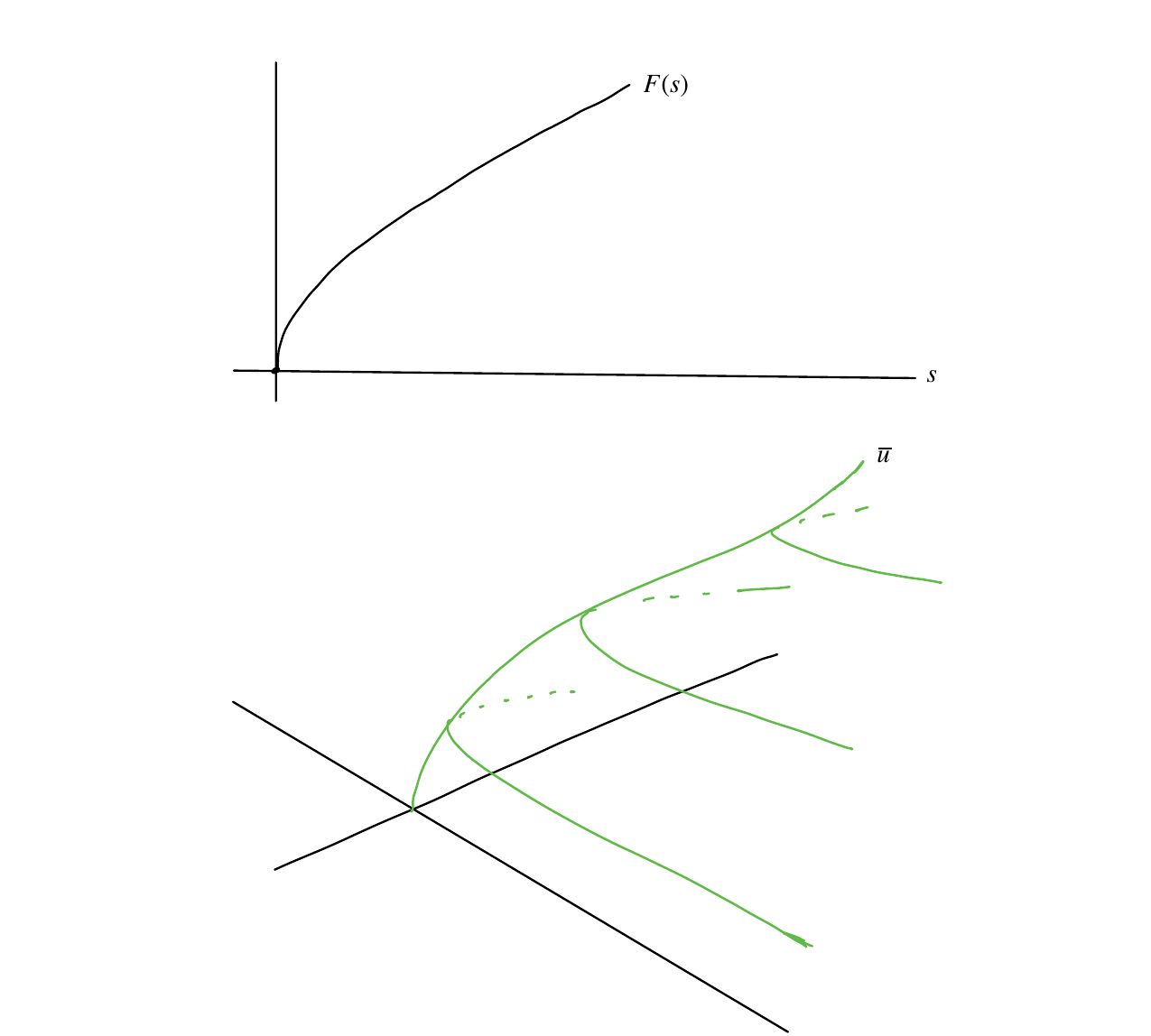}
\caption{}
\label{SuperSolnFig}
\end{figure}

\begin{rem}\label{SubSoln}
The function $\underline{u}$ can be obtained as follows. First, let $\underline{\sigma}(s) := 2(\sigma - \epsilon_0f)(s/2)$. Up to taking $\epsilon_0$ smaller (and increasing the constant $A$ accordingly in the construction of $\overline{u}$), we have $\overline{\sigma} < \underline{\sigma}$. Then replace $\overline{\sigma}$ by $\underline{\sigma}$ in the definitions of $w$ and $v$, replace $F'(s)$ by $G'(s) := e^{-B\int_s^{\infty} \frac{\tau^{-2/3}}{1+\tau^{2/3}}d\tau}$ for sufficiently large $B$, take $\underline{u} = \max\{G(v) - D,\,0\}$ in $\{|y| > |x|\}$ for sufficiently large $D$, and finally extend by odd reflection. The choices of $B$ and $D$ guarantee that $\underline{u}$ is a sub-solution to the MSE in $\{|y| > |x|\}$, and the ordering $\underline{u} \leq \overline{u}$ in $\{|y|>|x|\}$ follows easily from the ordering $\overline{\sigma} < \underline{\sigma}$ and the fact that $G'(s) < F'(s)$ for $s > 0$.
\end{rem}

\subsection{The Lawson-Osserman Cone}\label{LOEx}
In Section \ref{MSSBernstein} we noted that Lipschitz global solutions to the minimal surface system are linear provided one-homogeneous solutions are, and that the latter is true when the domain dimension is at most $3$. Here we briefly discuss the Lawson-Osserman example \cite{LO} of a four-dimensional graphical minimal cone of codimension three.

Recall that $\mathbb{S}^3 = \{z = (z_1,\,z_2): |z_1|^2 + |z_2|^2 = 1\} \subset \mathbb{C}^2 \cong \mathbb{R}^4$ can be identified with $SU(2)$ as follows:
$$z \leftrightarrow \left(\begin{array}{cc}
z_1 & z_2 \\
-\overline{z_2} & \overline{z_1}
\end{array}\right).$$
We let $\cdot$ denote the group operation that $\mathbb{S}^3$ inherits.
In turn, $SU(2)$ is isomorphic to the unit quaternions: if $z_1 = a + bi,\, z_2 = c+di$, take
$$q(z) = a + b\hat{i} + c\hat{j} + d\hat{k}.$$
Finally, there is a homomorphism from the unit quaternions to the rotations $SO(3)$ of $\mathbb{R}^3$. Identifying $\mathbb{R}^3$ with the pure quaternions $p$, this is given by $q \rightarrow R_q$, where
$$R_q(p) = qpq^{-1}.$$

The Hopf map $H$ from $\mathbb{S}^3$ to $\mathbb{S}^2$ is defined by
$$H(z) = R_{q(z)}\hat{i} = (|z_1|^2-|z_2|^2,\,-2iz_1z_2),$$
where we identify $\mathbb{R} \times \mathbb{C}$ with $\mathbb{R}^3$. We thus have, for any $w,\,z \in \mathbb{S}^3$,
$$(w \cdot z, R_{q(w)}H(z)) = (w \cdot z, H(w \cdot z)) \in \mathbb{S}^3 \times \mathbb{S}^2.$$
Thus, a portion of the graph of the map $u(z) = f(|z|)H(z/|z|)$ over any point $z \in \mathbb{C}^2 \backslash \{0\}$ can be taken by an isometry of $\mathbb{R}^7$, which splits as an isometry of the domain and target, into a portion of the same graph over $(|z|,\,0)$. In particular, to check whether $u$ solves the minimal surface system (\ref{MSS}), one only needs to check whether the system holds along $(t,\,0),\, t \in \mathbb{R}_+$. (In fact, the system becomes a nonlinear second-order ODE for $f$, see \cite{DY}).

In the particular case $f(|z|) = k|z|$, the graph is a cone, so to check minimality, one only needs to show that (\ref{MSS}) holds at the point $(1,0)$. It is easy to perform the calculations using the formula for $H$ above. Here are some details. At $(1,\,0)$ we have
$$g = \text{diag}(1+k^2,\,1,\,1+4k^2,\,1+4k^2).$$
The matrices $D^2u^{\alpha}(1,\,0)$ are zero along the diagonal for $\alpha = 2,\,3$, so (\ref{MSS}) holds at $(1,\,0)$ for $\alpha = 2,\,3$ and any $k$. We also have
$$D^2u^1(1,\,0) = k\text{diag}(0,\,1,\,-3,\,-3),$$
whence
$$g^{ij}u^1_{ij}(1,\,0) = k\left(1-\frac{6}{1+4k^2}\right).$$
This vanishes provided $k = \sqrt{5}/2$, i.e. the map
$$u(z) = \frac{\sqrt{5}}{2|z|}(|z_1|^2-|z_2|^2,\,-2iz_1z_2)$$
is a solution to the minimal surface system.

\begin{rem}
The Lawson-Osserman cone admits ``minimal de-singularizations" in the same spirit that the Simons cone does \cite{DY}. This involves a delicate study of the ODE for $f$ which guarantees that $u(z) = f(|z|)H(z/|z|)$ solves $(\ref{MSS})$.
\end{rem}

\newpage
\section{Recent Results and Further Directions}\label{RecentResults}

\subsection{Entire Minimal Graphs}
Simon \cite{Simon} showed that there are entire minimal graphs that blow down to cylinders over a variety of area-minimizing cones, including the area-minimizing Lawson cones. Nonetheless, many basic questions remain. For example, the known examples have polynomial growth. It would be interesting to decide whether all entire minimal graphs have polynomial growth. A related question is whether there exist nonlinear polynomial solutions. In the recent paper \cite{Guo}, polynomial solutions whose graphs blow down to the cylinders over all known examples of area-minimizing cones are ruled out, so new examples of algebraic area-minimizing cones would need to be constructed to answer this question in the positive.

\subsection{Anisotropic Minimal Hypersurfaces}
Another topic that has received attention recently is that of anisotropic minimal surfaces, namely, oriented hypersurfaces $\Sigma$ in $\mathbb{R}^{n+1}$ that are critical points of functionals of the form
\begin{equation}\label{Anisotropic}
A_{\Phi}(\Sigma) = \int_{\Sigma} \Phi(\nu)\,dA,
\end{equation}
where $\nu$ is the unit normal to $\Sigma$ and $\Phi$ is one-homogeneous, smooth and positive on $\mathbb{S}^n$, and has uniformly convex level sets. Important tools that are lost in this setting are the monotonicity formula (see \cite{All}) and rotation invariance.

Jenkins \cite{Jen} and Simon \cite{Simon2} proved that entire graphical minimizers of such functionals must be flat in dimensions $n = 2$ and $n = 3$, respectively. The anisotropic Bernstein problem was solved recently by constructing nonlinear entire graphical minimizers in the case $n = 4$ \cite{MY}, introducing the methods used in Section \ref{BDGEx} of this article. The analogue of the Simons cone in this construction is the cone over the Clifford torus in $\mathbb{R}^4$, which Morgan showed minimizes a functional of the type (\ref{Anisotropic}) \cite{Mor}. The examples in \cite{MY} have sub-quadratic growth. Polynomial entire anisotropic minimal graphs were constructed in the case $n = 6$ using a completely different method, based on solving a hyperbolic PDE \cite{M1}. It remains open whether there are polynomial solutions in lower dimensions. In the positive direction, flatness of entire anisotropic minimal graphs can be established in any dimension provided $\Phi$ is sufficiently close to the area integrand in an appropriate sense, and the gradient grows sufficiently slowly \cite{DYa}. Finally, the Morgan example shows that complete, stable critical points of such functionals need not be flat when $n \geq 3$. The question whether complete, stable critical points are flat remains open in the case $n = 2$ (see next sub-section for a discussion of this problem in the area case $\Phi|_{\mathbb{S}^n} = 1$).

\subsection{Stable Bernstein Problem}
An interesting problem in the theory of minimal hypersurfaces is to step away from graphicality and ask whether any complete, stable, two-sided minimal hypersurface in $\mathbb{R}^{n+1}$ is flat. A positive answer when $n = 2$ was given by Fischer-Colbrie and Schoen \cite{FCS}, Do Carmo and Peng \cite{DP}, and Pogorelov \cite{Pog2}. Schoen-Simon-Yau \cite{SSY} extended this up to dimension $n = 5$, and Bellettini up to dimension $n = 6$ \cite{Bel}, under the additional assumption that the volume growth is Euclidean (as it is e.g. for minimizers of area). The Simons cone says that $n \leq 6$ is necessary. Finally, a recent flurry of activity (\cite{CL2}, \cite{CMR}, \cite{CL1}, \cite{CLMS}, \cite{Maz}) has given a positive answer up to dimension $n = 5$, leaving only the case $n = 6$ unanswered.

\subsection{Complex Monge-Amp\`{e}re Equation}
As remarked above, there are non-quadratic global solutions to the complex Monge-Amp\`{e}re equation $\det \partial \bar{\partial} u = 1$ on $\mathbb{C}^n,\, n \geq 2$. It is natural to conjecture that global solutions with quadratic growth must be quadratic. This was established in \cite{MS} for solutions to the model equation $u_{11}u_{22} = 1$ on $\mathbb{R}^2$, which captures some of the structural features of complex Monge-Amp\`{e}re that present challenges (e.g. solutions are not convex, lack of invariance under certain rotations, invariance under adding certain quadratic polynomials). Here we are always assuming that $u$ satisfies the appropriate convexity condition so that the equation under study is elliptic, namely, plurisubharmonic for complex Monge-Amp\`{e}re, and convex in coordinate directions for $u_{11}u_{22} = 1$.


\end{document}